\documentclass[a4paper, 11pt, oneside]{article}

\usepackage{amsmath}
\usepackage[runin]{abstract}
\usepackage{hyperref}
\usepackage{fancyvrb}
\usepackage{amssymb}
\usepackage{mathtools}
\usepackage{amsthm}
\usepackage{amsfonts}
\usepackage{dsfont}
\usepackage{fancybox}
\usepackage{xcolor}
\usepackage{tikz}

\let\svthefootnote\thefootnote
\newcommand\freefootnote[1]{%
	\let\thefootnote\relax%
	\footnotetext{#1}%
	\let\thefootnote\svthefootnote%
}

\newtheorem*{algorithm*}{Algorithm}

\theoremstyle{plain}
\newtheorem{Theorem}{Theorem}[section]
\newtheorem{Lemma}[Theorem]{Lemma}
\newtheorem{Proposition}[Theorem]{Proposition}
\newtheorem{Definition}[Theorem]{Definition}

\theoremstyle{definition}
\newtheorem{Remark}[Theorem]{Remark}
\newtheorem{example}[Theorem]{Example} 

\numberwithin{equation}{section}

\newcommand{\supp}{{\rm supp \, }}

\author{Stephan Dahlke$^{\ast}$, Marc Hovemann$^{\ast}$, Thorsten Raasch$^{\dagger}$ and Dorian Vogel$^{\ast}$}
\title{Adaptive Quarklet Tree Approximation}
\date{\today}

\begin{document}
\maketitle

\freefootnote{
	\hspace{-21.5pt}
	This work has been supported by Deutsche Forschungsgemeinschaft (DFG), grants DA360/24 - 1 and RA2090/3 - 1.

\noindent	
${}^\ast$Philipps-Universit{\"a}t Marburg, Fachbereich Mathematik und Informatik,   Hans-Meerwein Str. 6,  Lahnberge, 35043 Marburg, Germany, Email: \texttt{dahlke@mathematik.uni-marburg.de}, \texttt{hovemann@mathematik.uni-marburg.de},
\texttt{vogeldor@mathematik.uni-marburg.de}.\\
\noindent	
${}^\dagger$Universit{\"a}t Siegen, Department Mathematik,  Walter-Flex-Str.
	3, 57068 Siegen, Germany, Email: \texttt{raasch@mathematik.uni-siegen.de}.	
	}
\noindent
\textbf{Abstract.} This paper is concerned with near-optimal approximation of a given function  $f \in L_2([0,1])$ with elements of a polynomially enriched wavelet frame, a so-called quarklet frame. Inspired by $hp$-approximation techniques of Binev, we use the underlying  tree structure of the frame elements to derive an adaptive algorithm that, under standard assumptions concerning the local errors, can be used to create approximations with an error close to the best tree approximation error for a given cardinality. We support our findings by numerical experiments demonstrating that this approach can be used to achieve inverse-exponential convergence rates.

\vspace{0,3 cm}

\vspace{0,3 cm}

\section{Introduction}

Many problems in science, engineering and finance are modeled as partial differential equations. Very often, a closed form of the unknown solution is not available, so that numerical schemes for its constructive approximation are indispensable tools. The most popular approach is the finite element method (FEM). The classical $h$-FEM relies on a space refinement of the domain under consideration. Alternatively, one can increase the polynomial degree of the ansatz functions. This is known as the $p$-method. A combination of both, so-called $hp$-FEM, is also possible. To handle large-scale real-life problems it is essential to employ \emph{adaptive} strategies increasing the overall efficiency. Then the goal is to end up with a satisfying approximation in a reasonable amount of time. In principle, an adaptive algorithm is an iterative strategy that identifies regions where the current approximation is still far away from the exact solution and refinement (polynomial enrichment) is performed only in these regions. In particular for adaptive $h$-FEM there is a huge amount of literature, we refer to \cite{bib:Cia02,bib:Hac10,NSV09,bib:Sch98}. Nowadays the convergence analysis of adaptive $p$- and $hp$-methods is more and more in the focus of research. These schemes converge very fast, often even exponentially. However, when it comes to the theoretical analysis and to rigorous convergence proofs with or without rates only a few results have arisen recently. To the best of our knowledge, the state of the art results concerning the convergence of $hp$-adaptive strategies are \cite{bib:BPS13,bib:BD11,bib:DV20,bib:DH07} and \cite{bib:CNSV14,bib:CNSV17}, which include optimality results. 

Another approach is to use wavelets. The advantages of wavelets are their strong analytical properties that can be used to derive adaptive schemes that are guaranteed to converge with the optimal convergence order of the best $N$-term wavelet approximation \cite{bib:CDD01,bib:Ste09}. Whenever the exact solution to the problem under consideration can be approximated by a linear combination of $N$ wavelets with an error proportional to $N^{-s}$,  where $s>0$ is the approximation rate, the adaptive wavelet scheme will be able to realize this very rate $s$. These schemes are essentially space refinement methods and can therefore be classified as $h$-methods. Then a very natural question is whether $hp$-versions of adaptive wavelet methods exist. At this point the concept of \emph{quarklets} comes into play. These polynomially enriched wavelets have been introduced in the last decade, see \cite{bib:DKR17}. They are constructed out of biorthogonal compactly supported Cohen-Daubechies-Feauveau spline wavelets, where the primal generator is a cardinal B-spline. For the theory of such biorthogonal wavelets we refer to Section 6.A in \cite{bib:CDF92}. Roughly speaking the quarklets are a linear combination of translated cardinal B-splines that are multiplied with some monomial. The precise definition can be found in Definition \ref{def_quarklet}. The properties of the quarklets have been studied in \cite{bib:DFK18,bib:DRS19,HoDa2021,bib:DHK22}. In particular, they can be used to design schemes that resemble $hp$-versions of adaptive wavelet methods. Furthermore it was shown in \cite{bib:DRS19} that it is possible to use quarklets to directly approximate certain model singularities that arise in elliptic boundary value problems on polygonal domains with the error exponentially decaying like $e^{-\alpha N^\beta}$, for some $\alpha,\beta>0$. It is our long term goal to design an adaptive scheme for solving partial differential equations based on quarklets that can be proven to converge and realizes the optimal (exponential) rate. This paper can be seen as a further step in this direction. 

Most adaptive schemes that strive for optimality share a common feature. In order to ensure optimal complexity they need to intertwine steps of refinement and derefinement. An example demonstrating this fact is given in Section 1.1 in \cite{bib:CNSV14}. Such a derefinement procedure is sometimes called \emph{coarsening}. The main idea can be summarized as follows. We consider a function, for example this could be the current approximation to the unknown solution of a partial differential equation, given as a linear combination of certain elements from our ansatz system. Now we approximate this function up to some tolerance with a rate close to the optimal convergence order. These coarsened approximations are generally less accurate than the original function but they possess an optimal balance between accuracy and degrees of freedom that can then be used to guarantee optimality of the overall scheme. For the application in $hp$-FEM a routine providing suitable approximations has recently been developed by Binev in \cite{Bin18}. There the $hp$-adaptive approximation is formulated as an approximation problem on a tree where the set of \emph{leaves} plays the role of a locally refined mesh and each leaf has a polynomial degree $p \ge 1$ assigned to it. In this paper we will design a similar routine in the quarklet setting. Therefore one needs certain structural demands on which quarklets can be used simultaneously for successful theoretical analysis and the practical proceeding. This constraint again arises in the form of a \emph{tree structure}. However, there is a major difference to classical $hp$-meshes. In the quarklet case the ancestors of the leaves, so-called \emph{inner nodes}, remain as an active contributor to the approximation. Consequently we have to assign a polynomial degree not only to the set of leaves but also to all inner nodes. The question how this is handled correctly is non-trivial. 

In this paper we will tackle this task and introduce a concept of \emph{quarklet trees} that works in theory and in practice. We use the theory developed in \cite{Bin18} as a starting point to obtain an algorithm {\bf NEARBEST\textunderscore TREE} that for a given function $f \in L_2([0,1])$ produces quarklet trees $T_N$ with cardinality $\# T_N \lesssim N$ that are proven to be near-best. By that we mean that a linear combination of the quarklets in the tree $T_N$ provides an approximation to $f$ with a global error close to the error of the optimal tree with cardinality $n \le N$, see Theorem \ref{theorem:1} for the main result. Moreover our numerical experiments show that for some natural univariate test cases we can indeed achieve exponential convergence rates with this routine. The authors are confident that {\bf NEARBEST\textunderscore TREE} can be used as a building block in an adaptive quarklet scheme for solving partial differential equations that can be proven to converge with a rate close to the optimal (exponential) rate.

This paper is organized as follows. In Section \ref{sec_quark} we recall the basic idea of quarklet frames as polynomially enriched wavelet bases and introduce a corresponding tree structure. Then, in Section \ref{sec:adap_ref} we show that the construction in \cite{Bin18} can be adapted to fit into the setting of quarklet frames. In particular, we can construct an algorithm that produces near-best tree approximations.  Finally, in Section \ref{sec_practice} we present one way to apply our algorithm in practice and conduct some first numerical experiments. 

Let us complete the introduction by fixing some notation. With  $c, C, C_1,C_2, \ldots $ we denote  positive constants. Unless stated otherwise they depend only on the fixed parameters. When we write $a\sim b$ we mean that there are constants $0 < C_1 \le C_2 < \infty$ such that $a \le C_{1} b \le C_2 a$. 

\section{Quarks and Quarklets} \label{sec_quark}

\subsection{Quarklets on the Real Line}\label{subsec_Quarklets_def1}
In this section we define quarks and quarklets. For that purpose in a first step we recall the definition of cardinal B-splines. Cardinal B-splines are defined by $ N_{1} \coloneqq \chi_{[0,1)}$ and for $ m \in \mathbb{N}  $ with $ m \geq 2  $ inductively by using the convolution
\begin{align*}
N_{m} \coloneqq N_{m - 1} \ast N_{1} = \int_{0}^{1} N_{m - 1}(\cdot - t) dt .
\end{align*}
In what follows for fixed $ m \in \mathbb{N} $ we will work with the symmetrized cardinal B-spline $ \varphi(x)  \coloneqq   N_{m}  ( x + \lfloor \frac{m}{2} \rfloor   )    $. We observe $ \supp \varphi = [ - \lfloor \frac{m}{2} \rfloor    ,  \lceil \frac{m}{2} \rceil    ]   $. The symmetrized cardinal B-spline shows up in the definition of the so-called quarks.
\begin{Definition}\label{Bquark}
	Let $ m \in \mathbb{N} $ and $ p \in \mathbb{N}_{0}  $. Then the \emph{$p$-th cardinal B-spline quark} $ \varphi_{p}  $ is defined by
	\begin{equation*}
	\varphi_{p}(x)  \coloneqq \Big ( \frac{x}{\lceil \frac{m}{2} \rceil } \Big )^{p}  N_{m} \Big ( x + \lfloor \frac{m}{2} \rfloor  \Big ) .
	\end{equation*}
\end{Definition}
The quarks will be very important for us in order to define the quarklets. Their properties have been studied in \cite{bib:DKR17}. For a given $ \tilde{m} \in \mathbb{N}  $ with $ \tilde{m} \geq m    $ and $  m + \tilde{m} \in 2 \mathbb{N}   $ there exists a compactly supported spline wavelet $ \psi  $ with
\begin{equation}\label{def_CDF_wav}
\psi = \sum_{k \in \mathbb{Z}} b_{k}  \varphi ( 2 \cdot - k   )
\end{equation}
with expansion coefficients $ b_{k} \in \mathbb{R}   $. Only finitely many of them are not zero. Moreover $ \psi  $ has $ \tilde{m}   $ vanishing moments and the  system
\begin{align*}
\Big \{ \varphi (  \cdot - k )  \ : \ k \in \mathbb{Z}  \Big  \} \cup \Big \{ 2^{\frac{j}{2}} \psi (2^{j} \cdot - k) \ : \ j \in \mathbb{N}_{0} \ , \ k \in \mathbb{Z}  \Big \}
\end{align*}
is a Riesz basis for $ L_{2}(\mathbb{R}) $. We refer to \cite{bib:CDF92} for details on the construction of $\psi$, see especially Section 6.A. We use these \emph{Cohen-Daubechies-Feauveau (CDF) spline wavelets} $ \psi $ to define the quarklets.  
\begin{Definition}\label{def_quarklet}
	Let $ p \in \mathbb{N}_{0}  $. Then the \emph{$p$-th quarklet} $ \psi_{p} $ is defined by 
	\begin{equation*}
	\psi_{p} \coloneqq \sum_{k \in \mathbb{Z}} b_{k} \varphi_{p}(2 \cdot - k).
	\end{equation*}
	Here the $ b_{k}  $ are the same as in \eqref{def_CDF_wav}. Moreover for $ j \in \mathbb{N}_{0}   $ and $ k \in \mathbb{Z}  $ we write 
	\begin{equation*}
	\psi_{p,j,k} \coloneqq 2^{\frac{j}{2}} \psi_{p}(2^{j} \cdot - k) \qquad \qquad \mbox{and} \qquad \qquad \psi_{p,-1,k} \coloneqq \varphi_{p}( \cdot - k) .
	\end{equation*}
\end{Definition}
It can be shown that the quarks and quarklets inherit some important properties of the B-splines and B-spline wavelets, respectively. In particular, Jackson and Bernstein estimates can be proved and the quarklets possess the same number of vanishing moments. For details concerning this topic we refer to \cite{bib:DKR17}. Suitably weighted quarklet systems form stable representation systems (so-called frames) for function spaces like the Lebesgue space $L_2(\mathbb{R})$ or the Sobolev spaces $H^s(\mathbb{R})$ with $0 \le s < m - \frac12$. Even more sophisticated spaces like the Besov or Triebel-Lizorkin spaces can also be characterized in terms of quarklets, see \cite{bib:sieber2020adaptive} and \cite{HoDa2021}. For readers' convenience we will recall the basic definition of a frame and state an alternative characterization. For a detailed overview on the topic of frames we refer to \cite{bib:Chr03}.
\begin{Definition}\label{def:frame}
    Let $\Lambda$ be a countable index set and $H$ be a Hilbert space with inner product $\langle \cdot, \cdot \rangle_H$. A system $\mathcal{F}=\{f_\lambda\}_{\lambda \in \Lambda} \subset H$ is called a  \emph{(Hilbert) frame} for $H$ if there exist constants $A,B>0$ such that for all $f \in H$ it holds
	\begin{equation*}\label{eq:frame}
		A \Vert f \Vert_{H}^2 \le \Vert\{ \langle f, f_\lambda \rangle_{H} \}_{\lambda \in \Lambda}\Vert_{\ell_2(\Lambda)}^2 \le B \Vert f \Vert_{H}^2.
	\end{equation*}
\end{Definition}
\begin{Proposition}\label{Prop:charc}
	A system $\mathcal{F} = \{f_\lambda\}_{\lambda \in \Lambda} \subset H$ is a frame for $H$ if and only if $\operatorname{clos}_H(\operatorname{span}(\mathcal{F})) = H$ and for all $f \in H$ it holds
	\begin{equation}\label{eq:frame_charac}
		B^{-1}\Vert f \Vert_H^2 \le \inf_{\{\mathbf{c} \in \ell_2(\Lambda) : f = \sum\limits_{\lambda \in \Lambda} c_\lambda f_\lambda \}} \Vert \mathbf{c} \Vert_{\ell_2(\Lambda)}^2 \le 	A^{-1}\Vert f \Vert_H^2.
	\end{equation}
\end{Proposition}
The proof of Proposition \ref{Prop:charc} can be found in \cite{Werner09}, see Proposition 2.2. In contrast to Riesz bases, frames allow for redundancy. Here we state the frame property of the quarklets in $L_2(\mathbb{R})$, see Theorem 3 in \cite{bib:DKR17}.
\begin{Theorem} \label{L_2_frame}
	Let the weights $w_p \ge 0$ be chosen such that $w_0=1$ and $w_p\left( p+1\right)^{-1/2}$ is summable. Then, the system
	\begin{equation*}
	\Psi_{Q,w} \coloneqq \{ w_p \psi_{p,j,k}: p \in \mathbb{N}_0, j \in \mathbb{N}_0 \cup \{-1 \}, k \in \mathbb{Z} \},
	\end{equation*}
	forms a frame for $L_2\left( \mathbb{R} \right) $.
\end{Theorem}

\begin{Remark}
Quarklet frames can also be constructed on quite general domains contained in $\mathbb{R}^d$. This is possible when working with boundary adapted quarklets. For this we refer to \cite{bib:DFK18}.
\end{Remark}

\subsection{Tree Structured Index Sets} \label{sec:TreeIndices}
In this section, we want to introduce tree structured index sets. Therefore in a first step we will recall the concept of wavelet indices. By $\Psi_0$ we denote the set of all wavelets
\begin{equation*}
\Psi_{0} \coloneqq \{ \psi_{0,j,k} :  j \in \mathbb{N}_{0}\cup\{-1\} , k \in \mathbb{Z}    \} 
\end{equation*}
and let $\Lambda_0$ be the set of all wavelet indices 
\begin{equation*}
    \Lambda_0 \coloneqq \{ (j,k) :  j \in \mathbb{N}_{0}\cup\{-1\} , k \in \mathbb{Z}    \} .
\end{equation*}
Consequently, a pair $ \lambda = (j,k) \in \Lambda_0  $ is called a \emph{wavelet index}. Those indices can be matched with reference intervals. For that purpose we define the dyadic intervals $ I_{j,k} \coloneqq 2^{-j}[k,k+1)$. Obviously each interval $ I_{j,k} $ refers to a wavelet index $(j,k)$ and vice versa. The intervals $ I_{j,k}  $ have some special properties that will be very important for us later: 
\begin{itemize}
	\item[(i)] For all $ (j,k) \in \Lambda_0 $ there exists a constant $ c > 0  $ independent of $j$ and $k$ such that we have $  \operatorname{supp} \psi_{0,j,k}  \subset  c  I_{j,k} $. In other words each interval $  I_{j,k}  $ and therefore each index $(j,k)$ can be associated with a wavelet $  \psi_{0,j,k}  $ and vice versa. 
	\item[(ii)] For fixed $ j \in \mathbb{N}_{0}  $ the intervals $  I_{j,k} $ form a disjoint partition of $ [0,1) $, namely $ \bigcup_{k = 0}^{2^{j}-1} I_{j,k} = [0,1)   $. Moreover let $ j_{1} \in \mathbb{Z}$ and $ k_{1} \in \mathbb{Z}   $ be fixed. Let $  j_{2} \in \mathbb{Z}   $ with $ j_{2} > j_{1}   $. Then there exist $ 2^{j_{2} - j_{1}}  $ intervals $ I_{j_{2},k}  $ such that there is the disjoint partition $ \bigcup_{k} I_{j_{2},k} = I_{j_{1},k_{1}}   $.
	\item[(iii)] The intervals $ I_{j,k}  $  are nested. Each interval $ I_{j,k} $ is a disjoint union of two intervals of the next higher level, i.e.,
    \begin{align*}
        I_{j,k} = I_{j+1 , 2k} \cup I_{j+1 , 2k+1} \qquad \qquad \mbox{and} \qquad \qquad I_{j+1 , 2k} \cap I_{j+1 , 2k+1} = \emptyset .
    \end{align*}	
    We will call $ I_{j+1 , 2k} $ and $ I_{j+1 , 2k+1} $ \emph{children} of $  I_{j,k} $. Conversely $  I_{j,k} $ is the \emph{parent} of $ I_{j+1 , 2k} $ and $ I_{j+1 , 2k+1} $.
    \item[(iv)] A combination of (ii) and (iii) yields that there also exist disjoint partitions of $  I_{j,k} $ using different $ j_{i} > j  $. In this context it is also possible to generate sequences of partitions of $  I_{j,k} $ by splitting one of the intervals from the partition into two intervals of the next finer level in each step. The partitions we obtain in this way are nested.    
\end{itemize}
Notice that each index $(j,k)$ has exactly two children (indices), namely $(j+1,2k)$ and $(j+1,2k+1)$. Therefore the concept of working with the intervals $ I_{j,k}  $ induces a natural ancestor-descendant relation. Hence, if for indices $ \lambda = (j_{1},k_{1}) \in \Lambda_0 $ and $ \mu = (j_{2},k_{2}) \in \Lambda_0 $ we have $ I_{j_{1},k_{1}} \subsetneq I_{j_{2},k_{2}}  $, we shall use the notation
\begin{equation*}
\lambda \succ \mu
\end{equation*}
and say that $\lambda$ is a \emph{descendant} of $\mu$. Conversely we will call $\mu$ an \emph{ancestor} of $\lambda$. By $\lambda \succeq \mu$ we mean that $\lambda$ is a descendant of $ \mu  $ or equal to $\mu$. Now we have the necessary tools to establish tree structured index sets.  

\begin{figure}[!ht]
	\centering
	\begin{minipage}[b]{0.49\textwidth} \centering
		\includegraphics[width=\textwidth]{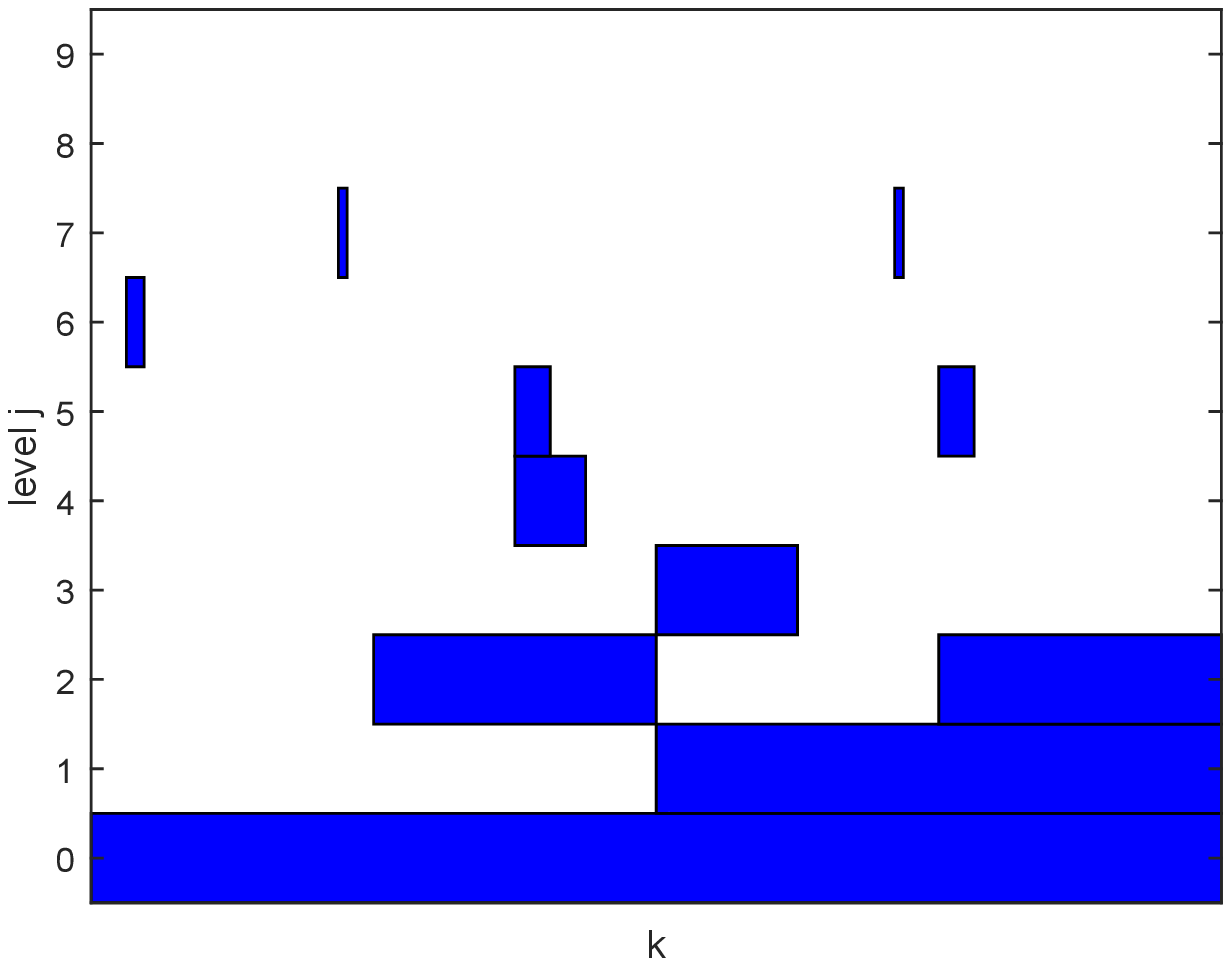}
		(a)
	\end{minipage} 
	\begin{minipage}[b]{0.49\textwidth} \centering
		\includegraphics[width=\textwidth]{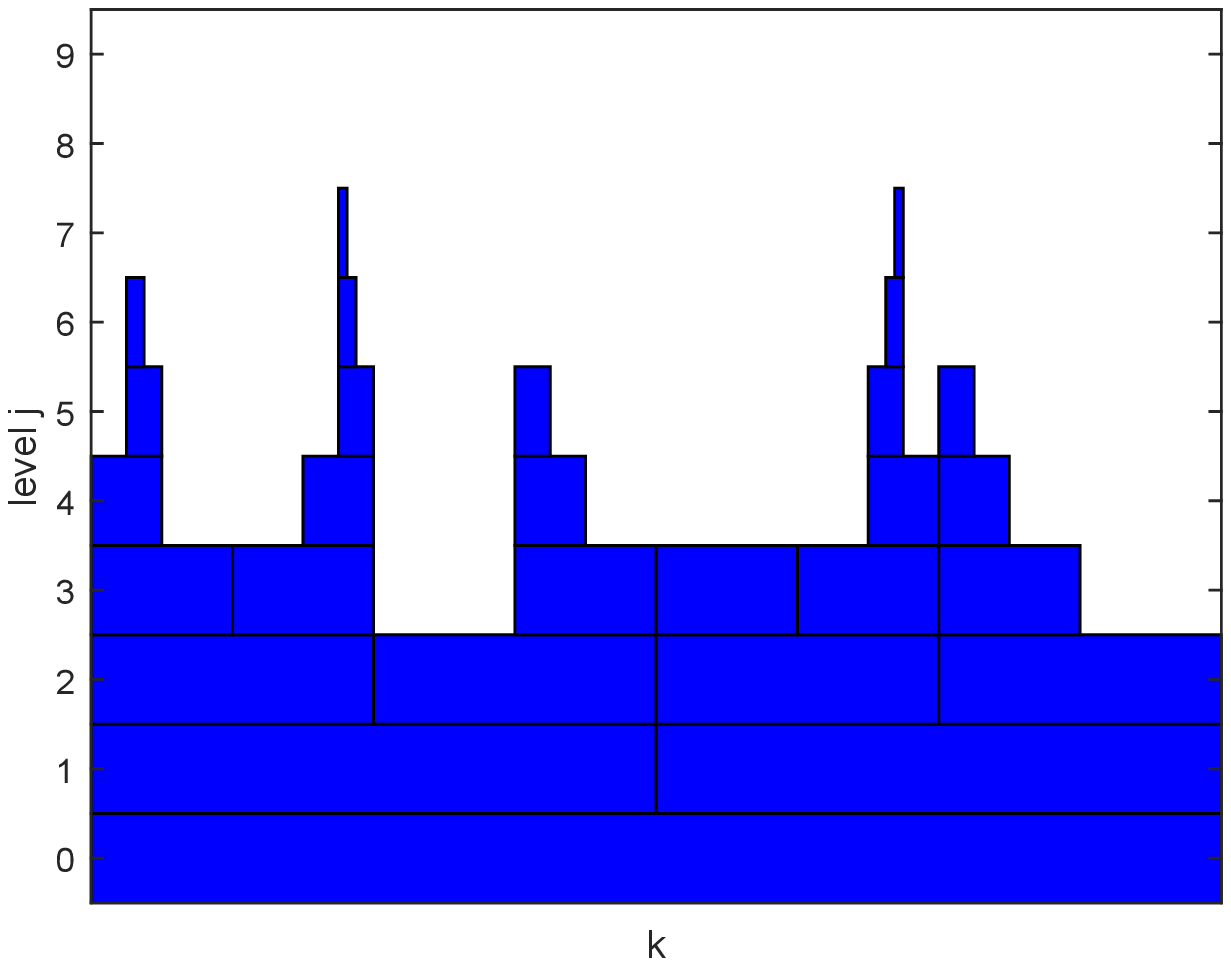}
		(b)
	\end{minipage}
    \begin{minipage}[b]{0.49\textwidth} \centering
	    \includegraphics[width=\textwidth]{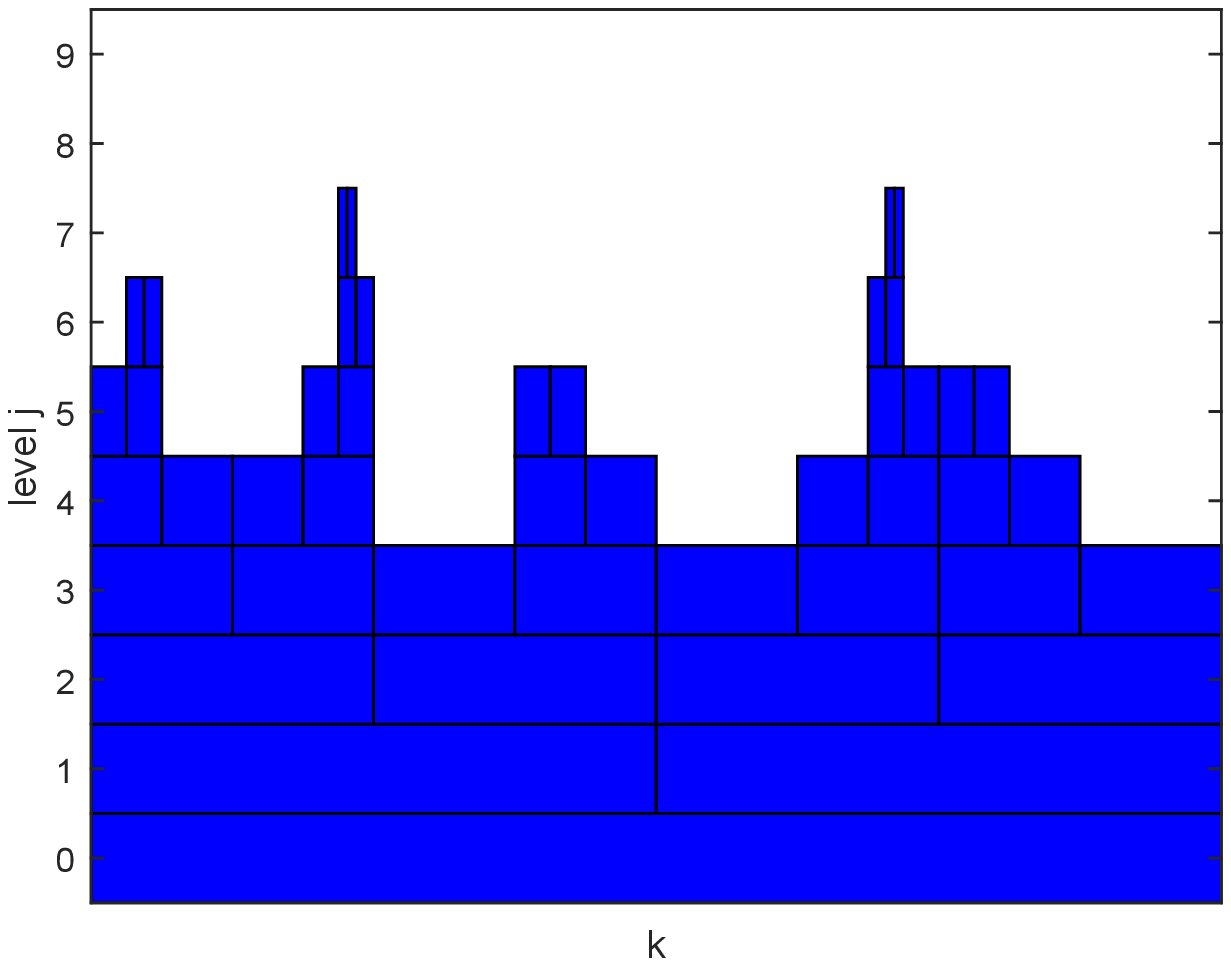}
	    (c)
    \end{minipage}
   \begin{minipage}[b]{0.49\textwidth} \centering
    	\includegraphics[width=\textwidth]{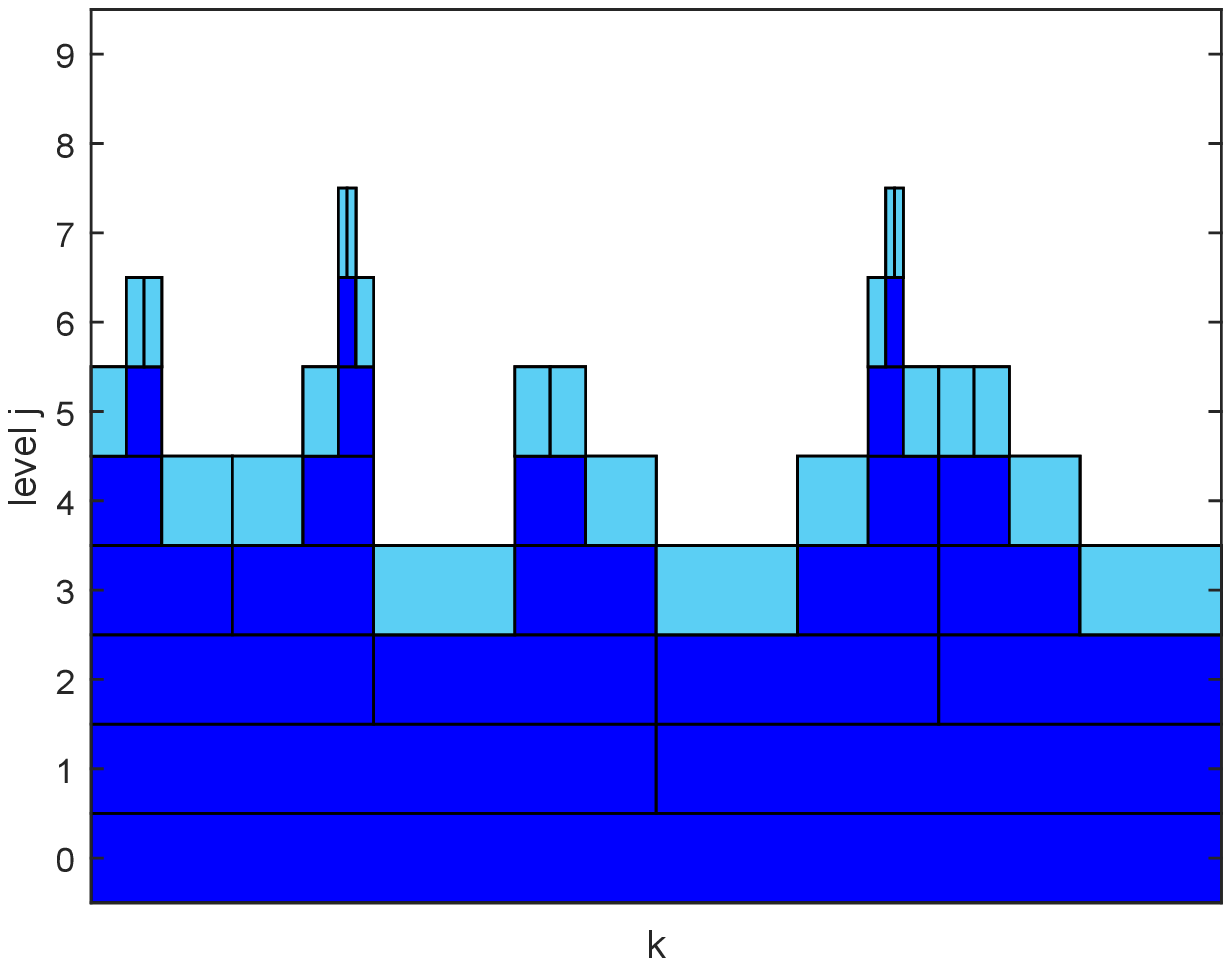}
    	(d)
   \end{minipage}
	\caption{An unstructured index set (a), the smallest (not complete) tree containing it (b), the corresponding completed tree (c) and the same tree again with highlighted leaves (d).}
	\label{Fig:index_sets}
\end{figure}

\begin{Definition}
Let $\mathcal{T} \subset \Lambda_0 $ be an index set. The set $  \mathcal{T} $  is called a \emph{tree (of wavelet indices)} if $\lambda \in \mathcal{T}$ implies $\mu \in \mathcal{T}$ for all $\mu \prec \lambda$. 
\end{Definition}
\begin{Definition} Let $\mathcal{T} \subset \Lambda_0   $ be a tree.
	\begin{itemize}
		\item[(i)]
		An index $\lambda \in \mathcal{T}$ is called \emph{node} of $\mathcal{T}$. 
		\item[(ii)]
		We set
		\begin{equation*}
		\mathcal{V}(\mathcal{T}) \coloneqq \{ \lambda \in \mathcal{T}  : \eta \notin \mathcal{T} \textrm{ for all } \eta \succ \lambda\}. 
		\end{equation*}
		The elements of $\mathcal{V}(\mathcal{T})$ are called \emph{leaves} of $  \mathcal{T} $. The set $\mathcal{T} \backslash \mathcal{V}(\mathcal{T})$ refers to the \emph{inner nodes}.
		\item[(iii)] An index $ \lambda \in \mathcal{T}   $ is called \emph{root} of $ \mathcal{T} $ if for all $ \eta \in \mathcal{T}   $ we have $\eta \succeq \lambda $. Then we often write $ \lambda = \mathcal{R}   $.
	\end{itemize}
\end{Definition}
Unless stated otherwise we will always work with \emph{complete} trees. That means an index $\lambda \in \mathcal{T}$ has exactly zero or two children, see Figure \ref{Fig:index_sets} for a visualization. For later use we also define $\mathcal{J}_{\lambda}$ to be the infinite wavelet tree rooted at $\lambda \in \Lambda_0$.

\begin{example}
The concept of wavelet indices we described above can be used to create a sequence of partitions of an interval that fits into a tree structure. For example let us consider $ I = [0,1]    $. Then the initial partition is the interval itself, namely $  I_{0,0} = [0,1]   $, and $  \mathcal{R} = (0,0)  $ is the root of the tree. To obtain the next partition we subdivide $  I_{0,0}  $ into two intervals of the next higher level, namely $ I_{1 , 0} $ and $ I_{1 , 1} $. This subdivision process can be repeated again and again to get a sequence of partitions. In each step we subdivide one element of the current partition, namely a leaf of the current tree, and attach two smaller intervals to it which are new leaves. The subdivided element is not a member of the new partition. But it remains in the tree as an inner node, which is not a leaf any more. In almost all steps of this process we have various choices which interval of the current partition we will subdivide to obtain the next partition. In our example to obtain the third partition we can subdivide either $ I_{1 , 0} $ or $ I_{1 , 1} $. In later steps the number of possible choices gets larger and larger. This will provide us the possibility to create an adaptive algorithm later. 
\end{example}

Now we are prepared to introduce local space refinement. For that purpose let $\mathcal{T} \subset \Lambda_0$ be a wavelet tree and $\lambda \in \mathcal{V}(\mathcal{T})$ be a leaf. The tree $\mathcal{T}$ is then refined by adding the two children $\eta_1, \eta_2 \in \Lambda_0$ of $\lambda$ to it. This space refinement can be written as
\begin{equation}\label{h-refinement}
\mathcal{T} \cup \{\eta_1, \eta_2\}.
\end{equation}
Of course this process can then be repeated with a different leaf.

\begin{Remark}
There are various ways to further generalize tree-structured index sets. In some cases they are even necessary for an efficient implementation. For example one can consider subdivisions with more than two children or start with multiple roots. Also adaptations for the multivariate case are possible \cite{bib:CDD03}.
\end{Remark}

Now we generalize the concept of wavelet indices to the more advanced concept of quarklet indices. To each wavelet index $  \lambda = (j,k) \in \Lambda_0 $ we match a third parameter $  p \in \mathbb{N}_{0}   $ to obtain a \emph{quarklet index} $(p,j,k) $. Each quarklet index refers to a quarklet $ \psi_{p,j,k}  $ and vice versa. By $\Psi$ we denote the set of all quarklets
\begin{equation*}
\Psi \coloneqq \{ \psi_{p,j,k} : p \in \mathbb{N}_{0} , j \in \mathbb{N}_{0}\cup\{-1\} , k \in  \mathbb{Z}    \} .
\end{equation*}
The corresponding index set is defined as
\begin{equation*}
\Lambda \coloneqq \{ (p,j,k) : p \in \mathbb{N}_{0} , j \in \mathbb{N}_{0}\cup\{-1\} , k \in  \mathbb{Z}     \} .
\end{equation*}
For a given quarklet index  $\lambda = (p,j,k) \in \Lambda$ we use the notation  $\vert \lambda \vert \coloneqq j $. Furthermore we use the mapping $\circ : \Lambda \rightarrow \Lambda_0 $ defined by
\begin{equation*}\label{eq:circ}
\lambda \mapsto \lambda^\circ \coloneqq (p,j,k)^\circ \coloneqq (0,j,k),
\end{equation*}
which provides the wavelet index for a given quarklet index. We will also consider wavelet indices as quarklet indices with $(j,k) = (0,j,k)$. Therefore we have $\Lambda_0 \subset \Lambda $. \\

Let us consider a wavelet tree $\mathcal{T} \subset \Lambda_0 \subset \Lambda$ as a set of quarklet indices. Now we have two options for the refinement of a leaf $\lambda \in \mathcal{V}(\mathcal{T})$. We can always refine in space by adding the two children wavelet indices of $\lambda$ to $\mathcal{T}$ as in \eqref{h-refinement}. Since we now work with quarklets we can also increase the polynomial degree by adding certain quarklet indices $(p,j,k)$ with $p>0$ to $\mathcal{T}$. However, the question which quarklets should be added is non-trivial, see Remark \ref{rem_p_ref_options} below. To this end we consider sets $\Upsilon(\lambda) \subset \Lambda_0$ such that for each wavelet tree $\mathcal{T} \subset \Lambda_0$ the union of the sets $\Upsilon(\lambda)$ over all leaves $\lambda \in \mathcal{V}(\mathcal{T})$ provides a disjoint decomposition of $\mathcal{T}$. More precisely we require the following properties on $\Upsilon(\lambda)$.
\begin{itemize}
	\item[(i)] The sets $\Upsilon(\lambda)$ have the form 
	\begin{equation}\label{eq:Upsilon}
		\Upsilon(\lambda) = \{\mu \in \Lambda_0 : \lambda \succeq \mu \succeq \mu_\lambda\}
	\end{equation} 
	for some fixed $\mu_\lambda \preceq \lambda$ independent of the tree $\mathcal{T}$.
    \item[(ii)] For each tree $\mathcal{T}$ and inner node $\lambda \in \mathcal{T}\backslash \mathcal{V}(\mathcal{T})$ with children $\eta_1,\eta_2 \in \mathcal{T}$ it holds $\Upsilon(\eta_1) \cap  \Upsilon(\eta_2) = \emptyset$.
    \item[(iii)] For each tree $\mathcal{T}$ it holds $\bigcup_{\lambda \in \mathcal{V}(\mathcal{T})} \Upsilon(\lambda) = \mathcal{T}$.
\end{itemize}
There are many ways to create the sets $\Upsilon(\lambda)$ depending on the precise rule for determining the node $\mu_\lambda$ in (i). We present one option in Section \ref{numerical_results}, see Figure \ref{Fig:hp-grid}. With the help of the sets $\Upsilon(\lambda)$ we can now introduce the polynomial enrichment of a leaf. We increase the maximal polynomial degree of $\lambda \in \mathcal{V}(\mathcal{T})$ by adding quarklet indices with $p = 1$ to each node $\mu \in \Upsilon(\lambda)$. More formally we can write this polynomial enrichment of the leaf $\lambda$ as
\begin{equation*}\label{eq:p-refinement}
\mathcal{T} \cup  \bigcup_{\mu = (i,\ell) \in \Upsilon(\lambda)} (p, i, \ell) .
\end{equation*}
This process can then be repeated with a different leaf or with the same leaf and $p=2$ and so forth.
\begin{Remark}\label{rem_p_ref_options}
	At a first glance the construction above does not seem to be the most intuitive way to incorporate polynomial enrichment. A straightforward idea might be to increase the polynomial degree only on the leaves of the tree under consideration like in \cite{Bin18}. However, this leads to sets of quarklets that often do not perform significantly better in practice than just the wavelet tree. To fully utilize the approximation power of quarklets with higher polynomial degrees it appears to be necessary to include quarklets with $p>0$ on different scales that maintain an ancestor-descendant relation. A very natural idea would then be to increase the maximal polynomial degree on each node in the path from the leaf $\lambda$ to the root. In this case the set $\Upsilon(\lambda)$ would simply be the set of all ancestors of $\lambda$ including $\lambda$ itself. But this means that the sets $\Upsilon(\lambda)$ are never disjoint because they all have at least the root as a common element. Unfortunately this in turn leads to severe problems for the theoretical analysis in Section \ref{sec:adap_ref}.
\end{Remark}
We are now ready to introduce quarklet trees.
\begin{Definition}\label{def:quarklet_tree}
	Let $T \subset \Lambda$ be an index set and for all $\lambda^\circ = (j,k) \in T$ define 
	\[ p_{\max}(\lambda^\circ) \coloneqq p_{\max}(\lambda^\circ,T) \coloneqq \max\{p \in \mathbb{N}_0 : (p,j,k) \in T\}. \]
	Then $  T $  is called a \emph{tree (of quarklet indices)} if the following conditions are fulfilled.
	\begin{itemize}
		\item[(i)]
		For each $\lambda \in T$ we have $\mu^\circ \in T$ for all $ \mu^\circ \preceq \lambda^\circ$. 
		\item[(ii)] For each $\lambda^\circ \in \mathcal{V}(T)$ we have $p_{\max}(\lambda^\circ) = p_{\max}(\mu^\circ)$ for all  $\mu^\circ \in \Upsilon(\lambda^\circ)$. 
		\item[(iii)] For each $\lambda^\circ = (j,k) \in T$ we have $(p,j,k) \in T$ for all $0 \le p \le p_{\max}(\lambda^\circ)$.
	\end{itemize}
\end{Definition}
In other words a quarklet tree structure consists of an underlying wavelet index set possessing a tree structure and moreover the nodes of this tree are enriched with all quarklets up to a certain polynomial degree. When we talk about a leaf, node or root of a quarklet tree we always mean the wavelet index, which is guaranteed to be an element of the tree and can be accessed from an arbitrary quarklet index via the mapping $\circ$.  By $|T|$ we denote the number of wavelet indices in an (arbitrary) index set $T$. For a quarklet tree $T$ we set its cardinality to be the number of quarklets in the tree, namely
\begin{equation*}
\#T \coloneqq |T| +  \sum_{\lambda^\circ \in T} p_{\textrm{max}}(\lambda^\circ).
\end{equation*}
We introduce another quantity with respect to a wavelet tree $\mathcal{T}$. For $\lambda \in \mathcal{T}$ we define
\begin{equation}\label{def_Q}
r(\lambda, \mathcal{T}) \coloneqq \vert \mathcal{V}(\mathcal{J}_{\lambda} \cap \mathcal{T}) \vert -1,
\end{equation}
which can be interpreted as the number of refinements that are necessary to create the subtree of $\mathcal{T}$ emanating from $\lambda$. If it is clear from the context we omit the second argument. This quantity will later be useful to undo certain steps of space refinement and instead employ $r(\lambda)$ steps of polynomial enrichment. \\

To finish this section we will have a look at two different ways to characterize a quarklet tree $T$. The first way is to consider a wavelet tree $\mathcal{T}$ and then fix the maximal polynomial degrees on all of its leaves. To this end we write
\begin{align*}
P_{\max} \coloneqq \big\{ p_{\textrm{max}}(\lambda)\big\}_{\lambda \in \mathcal{V}(\mathcal{T})}.
\end{align*}
Then the maximal polynomial degrees on the inner nodes are a direct consequence of Definition \ref{def:quarklet_tree}. For all $\lambda \in \mathcal{V}(\mathcal{T})$ and $\mu \in \Upsilon(\lambda)$ we have to set $p_{\max}(\mu) = p_{\max}(\lambda)$ to end up with a quarklet tree $T$. Therefore the assignment $P_{\max}$ already implies the maximal polynomial degree on all nodes not just on the leaves and we can write $T = (\mathcal{T}, P_{\max})$. For the second option we consider two wavelet trees $\mathcal{T}$, $\mathcal{T}'$ with $\mathcal{T} \subseteq \mathcal{T}'$. The subtree $\mathcal{T}$ implies a quarklet tree $T= (\mathcal{T},P_{\max})$ by setting $p_{\max}(\lambda) = r(\lambda,\mathcal{T}')$ for all $\lambda \in \mathcal{V}(\mathcal{T})$. Therefore we can also write $T = (\mathcal{T}, \mathcal{T}') = (\mathcal{T}, P_{\max})$ since $ P_{\max}$ is given by $\mathcal{T}$ and $\mathcal{T}'$. The intuition behind this is that we delete the descendants of $\lambda \in \mathcal{T}'$ and instead employ polynomial enrichment on $\lambda \in \mathcal{V}(\mathcal{T})$. This process is also called \emph{trimming} and will be useful for our adaptive scheme later on.

\subsection{Local Errors, Global Errors and Best Approximation}\label{subsec_localerr}
As stated above it is the main goal of this paper to construct an adaptive quarklet algorithm to approximate given functions in an efficient way. To describe what `efficient' means we have to introduce some error functionals. In a first step for each wavelet index $ \lambda = (j,k) \in \Lambda_0   $ and a polynomial degree $p \in \mathbb{N}_0$ we introduce local errors $ e_{p}(\lambda) : \Lambda_0 \rightarrow [0, \infty )    $. They are supposed to satisfy the following two very important properties.  
\begin{itemize}
	\item[(i)] There is a subadditivity of the error for the lowest order. That means for $\lambda \in \Lambda_0$ with children  $\eta_1$ and $\eta_2$ we require
	\begin{equation}\label{eq:1.2}
	e_0(\lambda) \ge e_0(\eta_1) +e_0(\eta_2).
	\end{equation}
	\item[(ii)] The error is reduced by increasing the polynomial degree,
	\begin{equation} \label{eq:1.3}
	e_p(\lambda) \ge e_{p+1}(\lambda). 
	\end{equation}
\end{itemize}
\begin{Remark}\label{rem_loc_err_f_seq}
The conditions \eqref{eq:1.2} and \eqref{eq:1.3}  are quite general. Nevertheless the local errors $e_p(\lambda)$ can be used to describe an adaptive quarklet algorithm in Subsection \ref{subsec_alg1} and to prove that this algorithm is near-best in Subsection \ref{subsec_alg_nearbest1}. However when it comes to practical applications, more conditions on the shape of the errors are required. Here one has different possibilities. Let $ I = [0,1]  $. Then on the one hand it is possible to define the local errors for functions $  f \in H^1(I)   $. An example of this can be found in Section 4.4 in \cite{bib:CNSV14}. The authors use $hp$-finite element methods to solve elliptic problems and expect the unknown solution to be in the Sobolev space $  H^{1}(I)$.  Another possibility is to define the local errors $ e_p(\lambda)  $ for quarklet coefficient sequences $ \mathbf{c} = \{ c_\lambda  \}_{\lambda \in \Lambda }  $ representing functions $   f \in L_2(I)  $, see Theorem \ref{L_2_frame}. This approach seems to be suitable for our purpose. More details and a precise definition of the local errors $  e_p(\lambda) \coloneqq e_p(\lambda , \mathbf{c})    $ for quarklet coefficient sequences can be found in Section \ref{sec:coeff_approach}.    
\end{Remark}

The local errors can be used to define global errors. For a given quarklet tree $T = (\mathcal{T}, P_{\max})$ we define the global error $\mathcal{E}(T)$ by
\begin{equation}\label{eq_glob_err}
\mathcal{E}(T) \coloneqq \sum_{\lambda \in \mathcal{V}(\mathcal{T})} e_{p_{\textrm{max}}(\lambda)}(\lambda).
\end{equation}
It collects the local errors for all leaves of the tree. One can also say that it picks up the local errors for all members of the current partition.

\begin{Remark}\label{rem_glob_err}
Let us remark that \eqref{eq_glob_err} can be understood in two different ways. On the one hand this equation just can be seen as a definition of the global error. On the other hand when it comes to practical applications one often needs that the global error has to satisfy additional requirements. For example, it should be equivalent to the norm of the residual. Such a demand has consequences for the choice of the local errors $  e_p(\lambda)  $. Therefore \eqref{eq_glob_err} also can be seen as an additional condition concerning the local errors. 
\end{Remark}

The global error can be used to define the so-called best approximation error. 
\begin{Definition}\label{def_bestappr_err}
The \emph{error of the best quarklet tree approximation} of cardinality $n \in \mathbb{N}$ is defined by
\begin{align*}
\sigma_n \coloneqq \inf_{T = (\mathcal{T}, P_{\max} )} \ \inf_{\# T \le n} \mathcal{E}(T).
\end{align*}
\end{Definition}
We can rewrite this best approximation error as
\begin{equation*}
\sigma_n = \inf_{T = (\mathcal{T},\mathcal{T}')} \ \inf_{ \#T \le n } \mathcal{E}(T).
\end{equation*}
This motivates the approach of finding a wavelet tree $\mathcal{T}'$ first and then examine all possible subtrees $\mathcal{T} \subseteq \mathcal{T}'$. It will be our main goal to find an incremental algorithm that for each $N \in \mathbb{N}$ produces a quarklet tree $  T_N = (\mathcal{T}_N,\mathcal{T}'_N) = (\mathcal{T}_N, P_{\max}) $ with $  \# T_N \le \tilde{C} N $  that provides a near-best quarklet approximation in the sense of
\begin{equation}\label{eq:nearbest}
\mathcal{E}(T_N) \le C \sigma_{c N},
\end{equation}
with independent constants $C \ge 1 $ and $c \in (0,1]$. In other words, we are looking for an algorithm that produces a quarklet tree approximation with a global error comparable to the best possible error. For a detailed discussion concerning the constants $c,C,\tilde{C}$ we refer to Lemma \ref{lemma:T_N_cardinality} and Remark \ref{remark:interpretation}.

\section{Adaptive Refinement Strategy} \label{sec:adap_ref}

\subsection{Error Functionals for Adaptive Refinement}
To construct our near-best quarklet algorithm we need some more error functionals as well as indicators to decide where to refine a given tree. We will split their introduction into three steps. In the first step we will have a look at a penalized version of the local error for the lowest order that can be used to design near-best space adaptive schemes. Then we will define a counterpart for the space and polynomial degree adaptive case. In the last step we introduce two indicators that will help us decide where to refine based on these information. All these functions are specially tailored for our purposes and play a key role in both the algorithm and the proof of its optimality. \\

\noindent
\textit{Step 1:} We will need a modified local error functional denoted by $\tilde{e}(\lambda)$ with $ \tilde{e}(\lambda) : \Lambda_0 \rightarrow [0, \infty )    $. It is strongly connected with the local error of the lowest order $ e_{0}(\lambda)  $ and will be defined recursively. For the sake of convenience we write $e(\lambda) = e_0(\lambda)$. Let $ \mathcal{R}  $ be the root of a tree and $\mu$ be the parent of $\lambda$. Then we define the modified errors $ \tilde{e}  $ step by step out of the local errors via
\begin{equation}\label{eq:2.3}
\tilde{e}(\mathcal{R})\coloneqq e(\mathcal{R}), \qquad \tilde{e}(\lambda) \coloneqq \frac{e(\lambda)\tilde{e}(\mu)}{e(\lambda)+\tilde{e}(\mu)}.
\end{equation}
In the case $e(\lambda)=\tilde{e}(\mu)=0$ we put $\tilde{e}(\lambda)=0$. An equivalent formula to (\ref{eq:2.3}) is given by
\begin{equation*}\label{eq:2.4}
\frac{1}{\tilde{e}(\lambda)} = \frac{1}{e(\lambda)}+\frac{1}{\tilde{e}(\mu)}.
\end{equation*}
Iterating this equation leads to the useful relation
\begin{equation}\label{eq:2.7}
\frac{1}{\tilde{e}(\lambda)} = \sum_{\mu \preceq \lambda} \frac{1}{e(\mu)}.
\end{equation}

\noindent
\textit{Step 2:} We introduce an error functional denoted by $ E(\lambda)   $ with $ E(\lambda) : \mathcal{T} \rightarrow [0 , \infty )  $. It is only defined for a given wavelet tree $ \mathcal{T}  $, which is why we sometimes write  $E(\lambda) \coloneqq E(\lambda, \mathcal{T})$. It is explained recursively starting at the leaves of a tree $\mathcal{T}$. Here for $ \lambda \in \mathcal{V}(\mathcal{T})  $ we define
\begin{equation*}\label{eq:3.0}
E(\lambda)\coloneqq e(\lambda) = e_0(\lambda).
\end{equation*}
For the inner nodes of the tree this error functional is defined step by step moving from the leaves towards the root. Let $\eta_1$ and $\eta_2$ be the children of $\lambda$ and assume that $E(\eta_1)$ and $E(\eta_2)$ are already defined. Moreover recall that  $ r(\lambda) = r(\lambda, \mathcal{T})   $ is given by \eqref{def_Q}. Then for an inner node $ \lambda  $ we put    
\begin{equation}\label{eq:3.1}
E(\lambda) \coloneqq \min\{E(\eta_1)+E(\eta_2), e_{r(\lambda)}(\lambda)\}.
\end{equation}
This refers to the adaptive choice between the two refinement types. 

Next we want to introduce a modified version of $E(\lambda)$. To this end we first notice that enlarging the tree $\mathcal{T}$ changes the quantity $E(\lambda) = E(\lambda, \mathcal{T})$ only if $r(\lambda, \mathcal{T})$ changes. A similar observation was made in \cite{Bin18}, see page 3353. Here also some more explanations concerning this topic can be found, see Remark 3.1. We use this observation and consider a sequence $\mathcal{T}_1,\mathcal{T}_2,\mathcal{T}_3,\ldots$ of growing trees. With that we mean that each tree $\mathcal{T}_{k+1}$ is derived from $\mathcal{T}_k$  by subdividing a leaf and adding the two children indices. For a node $\lambda$ and $j \in \mathbb{N}_0$ there might exist multiple trees $\mathcal{T}_\star$ with $\lambda \in \mathcal{T}_\star$ and $r(\lambda, \mathcal{T}_\star)=j$ in the sequence $\mathcal{T}_1,\mathcal{T}_2,\mathcal{T}_3,\ldots$ of trees. This means that the subtree emanating from $\lambda$ stays the same in all the trees $\mathcal{T}_\star$ and consequently the quantity $E(\lambda,\mathcal{T}_\star)$ does not change. By using this observation we can let $ j \in \mathbb{N}_0  $ and $\mathcal{T}_\star$ be any of the trees in the sequence $\mathcal{T}_1,\mathcal{T}_2,\mathcal{T}_3,\ldots$ such that $r(\lambda, \mathcal{T}_\star)=j$ to define
\begin{equation}\label{eq:3.1a}
E_j(\lambda) \coloneqq E(\lambda, \mathcal{T}_\star) . 
\end{equation}
Using the error functionals $  E_j(\lambda)  $ as a starting point, we can now define  modified errors  $\tilde{E}_j(\lambda)$. They have some similarities with those explained in formula \eqref{eq:2.3}. For $ j = 0  $ we put $\tilde{E}_0(\lambda)\coloneqq \tilde{e}(\lambda)$. For $ j > 0  $ the error functionals $ \tilde{E}_j(\lambda)  $ are defined recursively via
\begin{equation}\label{eq:3.2}
\tilde{E}_j(\lambda)\coloneqq \frac{E_j(\lambda)\tilde{E}_{j-1}(\lambda)}{E_j(\lambda)+\tilde{E}_{j-1}(\lambda)}.
\end{equation}
For the case $E_j(\lambda)=\tilde{E}_{j-1}(\lambda)=0$ we put $\tilde{E}_j(\lambda)=0$. The error functional $ \tilde{E}_j(\lambda)  $  can be rewritten in terms of some of the other error functionals we explained above. We use the definition of $\tilde{E}_j(\lambda)$ several times and plug in formula \eqref{eq:2.7} to observe
\begin{align}
\frac{1}{\tilde{E}_j(\lambda)} = \frac{1}{E_j(\lambda)}+\frac{1}{\tilde{E}_{j-1}(\lambda)} =\sum_{k=1}^j \frac{1}{E_k(\lambda)}+\frac{1}{\tilde{E}_0(\lambda)}
= \sum_{k=0}^j \frac{1}{E_k(\lambda)} + \sum_{\mu \prec \lambda} \frac{1}{e(\mu)}. \label{eq:3.3}
\end{align}
As an outgrowth of $ \tilde{E}_j(\lambda)   $ we use \eqref{eq:3.2} with $ j = r(\lambda)   $ to set 
\begin{equation*}
\tilde{E}(\lambda) \coloneqq \tilde{E}(\lambda, \mathcal{T}) \coloneqq \tilde{E}_{r(\lambda)}(\lambda).
\end{equation*}

\noindent
\textit{Step 3:} We will need two more functions denoted by $q$ and $s$. They are strongly connected with the error functionals we defined above. At first let us define the function $q: \mathcal{T} \rightarrow [0,\infty)$. For a leaf $\lambda \in \mathcal{V}(\mathcal{T})$ it is just the modified local error, namely
\begin{equation*}
q(\lambda) \coloneqq \tilde{e}(\lambda) = \tilde{E}_0(\lambda).
\end{equation*}
For an inner node $\lambda \in \mathcal{T} \backslash \mathcal{V}(\mathcal{T})$ the function $q$ is defined recursively moving from $ \lambda  $ towards the leaves. With $\eta_1$ and $\eta_2$ denoting the children of $ \lambda  $, we put
\begin{equation}\label{def_func_q1}
q(\lambda)  \coloneqq \min \left\{\max\{q(\eta_1),q(\eta_2)\}, \tilde{E}_{r(\lambda)}(\lambda)\right\}.
\end{equation}
Roughly speaking one can say that $ q(\lambda)  $ is connected with the error resulting out of the subtree $ \mathcal{T} \cap \mathcal{J}_{\lambda}    $. Now let us define the function $s:\mathcal{T} \rightarrow\mathcal{V}(\mathcal{T})$. It maps each node of the tree to a leaf of it. For a leaf $\lambda \in \mathcal{V}(\mathcal{T})$ itself we just set
\begin{equation*}
s(\lambda) \coloneqq \lambda.
\end{equation*} 
For an inner node $\lambda \in \mathcal{T} \backslash \mathcal{V}(\mathcal{T})$ the function $s$ is defined recursively by  using the function $q$. With $\eta_1$ and $\eta_2$ denoting the children of $ \lambda  $, we put
\begin{equation*}\label{eq:3.5}
s(\lambda) \coloneqq s\left( \operatorname{argmax}\{ q(\eta_1),q(\eta_2)\} \right) .
\end{equation*}
Hence, $s$ points to that leaf of the investigated subtree with the largest penalized local error. Later on $s$ plays a key role for our algorithm since it tells us which leaf will be refined next. 

\begin{Remark}
	All the error functionals $ \tilde{e}(\lambda)  $, $  E(\lambda)  $ and $ \tilde{E}(\lambda)   $ are based on the local errors $  e_p(\lambda)  $, see \eqref{eq:1.2} and \eqref{eq:1.3}. Therefore the way how we define $  e_p(\lambda)  $ has direct consequences for the shape of the modified error functionals.  
\end{Remark}

\begin{Remark}
	The modified local errors from \eqref{eq:2.3} can be used to develop space adaptive refinement strategies. Indeed, they have already been successfully employed in adaptive wavelet and finite element methods for operator equations, see \cite{Bin18}, \cite{binev2004adaptive}, \cite{binev2004fast} and \cite{kappei2011adaptive}. 
\end{Remark}

\subsection{An Algorithm for Quarklet Tree Approximation}\label{subsec_alg1}

Now we are ready to explicitly state our adaptive quarklet algorithm. It is based on the local errors $   e_p(\lambda)    $ and the modified local error functionals defined upon them. Without loss of generality, we consider functions on $I = [0,1]$, by suitable rescaling arguments, the general case can be reduced to this model setting. Then as an input we can either use a function $f \in L_2(I)$ or a sequence of its quarklet expansion coefficients. We use the notation $\mathbf{f}$ which  stands for either of the two options depending on the definition of $e_p(\lambda)$. The algorithm {\bf NEARBEST\textunderscore TREE} takes  $ \mathbf{f} $ as input and adaptively produces a tree $  \mathcal{T}'_N  $ consisting of wavelet indices only. This tree is designed in such a way that it can easily be transformed into a tree $T_N = ( \mathcal{T}_N, P_{\max}) = (\mathcal{T}_N, \mathcal{T}'_N) $ of quarklet indices that provides a near-best quarklet tree approximation. The algorithm itself is similar to that in \cite{Bin18}. Let us remark that we do not fix the root $\mathcal{R}$ here. However, in our applications we will always use $\mathcal{R} = (0,0)$.

\vspace{0,4 cm}
\medskip
\noindent\shadowbox{\parbox{0.96\textwidth}{
		\begin{algorithm*} {\bf NEARBEST\textunderscore TREE} $[\mathbf{f},N_{\max}] \mapsto \mathcal{T}'_N$\\[-3ex]
			\noindent{
				\begin{tabbing}
					set $\mathcal{T}'_0 \coloneqq \{\mathcal{R}\}$, $\tilde{e}(\mathcal{R}) \coloneqq e(\mathcal{R})$, $E_0(\mathcal{R}) \coloneqq e(\mathcal{R})$, $\tilde{E}_0(\mathcal{R}) \coloneqq \tilde{e}(\mathcal{R})$,\\[0.1cm]
					$q(\mathcal{R}) \coloneqq \tilde{e}(\mathcal{R})$, $s(\mathcal{R}) \coloneqq \mathcal{R}$, $r(\mathcal{R}) \coloneqq 0$;\\[0.1cm]
					\texttt{\em for}\={} $N  = 1$ \texttt{to} $N_{\max}$ \\[0.1cm]
					\>  expand the current tree $\mathcal{T}'_{N-1}$ to $\mathcal{T}'_{N}$ by subdividing $\lambda_N \coloneqq s(\mathcal{R})$ and \\[0.1cm]
					\> adding its children $\hat{\eta}_1$ and $\hat{\eta}_2$ to it; \\[0.1cm]
					\> \texttt{\em for}\={} $\lambda \in  \{\hat{\eta}_1, \hat{\eta}_2\}$ \\
					\> \> calculate $\tilde{e}(\lambda) \coloneqq \frac{e(\lambda)\tilde{e}(\lambda_N)}{e(\lambda) + \tilde{e} (\lambda_N)}$, $E_0(\lambda) \coloneqq e(\lambda)$,  $\tilde{E}_0(\lambda) \coloneqq \tilde{e}(\lambda)$,\\[0.1cm]    
					\> \> $q(\lambda) \coloneqq \tilde{e}(\lambda)$, $s(\lambda) \coloneqq \lambda$, $r(\lambda) \coloneqq 0$; \\[0.1cm]
					\> \texttt{\em end for} \\
					\> set $\lambda = \lambda_N$; \\[0.1cm]
					\> \texttt{\em whi}\=\texttt{\em le} $\lambda \neq \emptyset$ \\[0.1cm]
					\> \> set $r(\lambda) \coloneqq r(\lambda) + 1 $ and calculate $e_{r(\lambda)}(\lambda)$;\\[0.1cm]
					\> \> set $\eta_1$ and $\eta_2$ to be the children of $\lambda$; \\[0.1cm]
					\> \> set $E_{r(\lambda)}(\lambda) \coloneqq \min\{E_{r(\eta_1)}(\eta_1)+E_{r(\eta_2)}(\eta_2), e_{r(\lambda)}(\lambda)\}$; \\[0.1cm]
					\> \> set $\tilde{E}_{r(\lambda)}(\lambda)\coloneqq \frac{E_{r(\lambda)}(\lambda)\tilde{E}_{{r(\lambda)}-1}(\lambda)}{E_{r(\lambda)}(\lambda)+\tilde{E}_{{r(\lambda)}-1}(\lambda)}$; \\[0.1cm]
					\> \> set $\eta \coloneqq \operatorname{argmax}\{q(\eta_1), q(\eta_2)\}$, $q(\lambda) \coloneqq \min\{q(\eta), \tilde{E}_{r(\lambda)}(\lambda)\}$ and $s(\lambda) \coloneqq s(\eta)$; \\[0.1cm]
					\> \> replace $\lambda$ with its parent (or $\emptyset$ if $\lambda = \mathcal{R}$); \\[0.1cm]
					\> \texttt{\em end while}\\[0.1cm]
					\texttt{\em end for}
				\end{tabbing}
			}
		\end{algorithm*}
}}
\medskip
\vspace{0,4 cm}

The main idea of the algorithm {\bf NEARBEST\textunderscore TREE} can be summarized in the following way. We start with the tree $\mathcal{T}'_0 \coloneqq \{\mathcal{R}\}$ and then proceed iteratively. As long as we have $N \le N_{\textrm{max}}$ for the given tree $\mathcal{T}'_{N-1}$ we subdivide the leaf $s(\mathcal{R})$ and add two child nodes to it to form $\mathcal{T}'_{N}$. The inner loop updates all the important quantities going from the newly created leaves back to the root $\mathcal{R}$. Notice that on all nodes which are not on this path no changes are needed. At this point let us stress the significance of the error functional $\tilde{E}$. Since we are interested in a reduction of the error $ E(\mathcal{R})$, a straightforward greedy strategy based on the errors $E(\lambda)$ may fail to guarantee an effective error reduction. To overcome this problem we work with the modified error $\tilde{E}(\lambda)$. It penalizes the lack of success in the reduction of $E(\lambda)$, and then $s(\mathcal{R})$ points to the leaf with the largest penalized error. For that reason $\tilde{E}$ and the subsequent decisions made in $q$ and $s$ based on its information can be considered as the adaptive steering wheel in this algorithm. The following lemma provides an estimate of the complexity of the algorithm {\bf NEARBEST\textunderscore TREE}.

\begin{Lemma}\label{lem_near_best_complex}
Let $N \in \mathbb{N}$. The incremental algorithm {\bf NEARBEST\textunderscore TREE}	performs $\sum_{\lambda \in \mathcal{T}'_N} (r(\lambda)+1)$ steps to obtain $\mathcal{T}'_N$.
\end{Lemma}
This result can be shown by following the lines of the proof of Lemma 3.2 on page 3354 in \cite{Bin18}. Here a similar algorithm has been investigated for $hp$-adaptive approximation. Therefore we skip the details of the proof. The complexity identified in Lemma \ref{lem_near_best_complex} depends on the balancing of the tree. It varies from $\mathcal{O}(N \log N)$ to $\mathcal{O}(N^2)$ in the best and worst case, respectively. 

Now let us have a closer look at the process of trimming. As already mentioned the tree $ \mathcal{T}'_N  $ produced by the algorithm {\bf NEARBEST\textunderscore TREE} consists of wavelet indices only. However, it is designed such that we can transform it into a quarklet tree $T_{N} = (\mathcal{T}_N, \mathcal{T}'_N)$ easily. To obtain the optimal subtree $\mathcal{T}_N$ of $\mathcal{T}'_N$ we start with $\mathcal{T}'_N$. Then we move from the root $ \mathcal{R}$ toward a leaf $\eta \in \mathcal{V}(\mathcal{T}'_N)$. At the first node where we observe $E(\lambda) = e_{r(\lambda)}(\lambda)$ in \eqref{eq:3.1} we trim the tree. That means we delete all descendants of $\lambda$. By definition we have $E(\eta) = e_0(\eta)$ on the leaves $\eta \in \mathcal{V}(\mathcal{T}'_N)$. Consequently, we are guaranteed to encounter this situation. This procedure is then repeated for all remaining paths which have not been treated. This way $\mathcal{T}_N$ becomes the minimal tree for which we have $E(\lambda) = e_{r(\lambda)}(\lambda)$ on all its leaves. The tree and its subtree can be used to define a quarklet tree $T_N  = (\mathcal{T}_N, \mathcal{T}'_N) = (\mathcal{T}_N, P_{\max})$ by setting $p_{\max}(\lambda) = r(\lambda, \mathcal{T}'_N)$ as explained in Section \ref{sec:TreeIndices}. If $\mathcal{T}'_N$ stems from the algorithm {\bf NEARBEST\textunderscore TREE} the quantities $E(\lambda)$ and  $e_{r(\lambda)}(\lambda)$ can be directly extracted from the algorithm. Then no further calculations are needed. The following algorithm {\bf TRIM} provides one way to implement the trimming procedure.

\vspace{0,4 cm}

\medskip
\noindent\shadowbox{\parbox{0.96\textwidth}{
		\begin{algorithm*} {\bf TRIM} $[\mathcal{T}'] \mapsto \mathcal{T}$\\[-3ex]
			\noindent{
				\begin{tabbing}
					set $B = \{\mathcal{R}\}$ and $ \mathcal{T} = \mathcal{T}'$;\\[0.1cm]
					\texttt{\em whi}\=\texttt{\em le} $B \neq \emptyset$ \\[0.1cm]
					\> take $\lambda \in B$; \\[0.1cm]
					\> \texttt{\em if}\={} $E(\lambda) = e_{r(\lambda)}(\lambda)$\\[0.1cm]
					\> \> remove all descendants from $\lambda$ in $\mathcal{T}$;\\[0.1cm]
					\> \texttt{\em else}\\[0.1cm]
					\> \> add the children $\eta_1$ and $\eta_2$ of $\lambda$ to $B$;\\[0.1cm]
					\> \texttt{\em end if}\\[0.1cm]
					\> remove $\lambda$ from $B$;\\[0.1cm]
					\texttt{\em end while}
				\end{tabbing}
			}
		\end{algorithm*}
}}
\medskip

\vspace{0,4 cm}

Next, we estimate the cardinality of the tree $T_N$ created by our algorithm.
\begin{Lemma} \label{lemma:T_N_cardinality}
	Let $  N \in \mathbb{N}  $ and let $  \mathcal{R} = (0,0)$.	Let the tree $T_N = (\mathcal{T}_N,\mathcal{T}'_N)$ be produced by the algorithm {\bf NEARBEST\textunderscore TREE} and a subsequent trimming. Then it holds
	\begin{equation}\label{T_N_cardinality}
	N+1 \le \# T_N \le \frac{N^2+6N+5}{4}
	\end{equation}
	and
	\begin{equation}\label{T_N_cardinality2}
	\# T_N \le \left( \max_{\lambda \in \mathcal{T}_N} |\lambda|+1\right) N+1.
	\end{equation}
\end{Lemma}
\begin{proof}
	Let us start by recalling that each quarklet tree $T$ can be characterized by two sequences of refinements. The first one describes the refinements in space in form of a wavelet tree $\mathcal{T}$, while the second one expresses the steps of polynomial enrichment on $\mathcal{T}$ in terms of an assignment $\{ p_{\textrm{max}}(\lambda)\}_{\lambda \in \mathcal{V}(\mathcal{T})}$. Now let $N_h \coloneqq r(\mathcal{R},\mathcal{T})$ and $N_p \coloneqq \sum_{\lambda \in \mathcal{V}(\mathcal{T})} p_{\max}(\lambda)$ with $N=N_h+N_p$ denote the total numbers of refinements in space and polynomial degree, respectively. Adding two finer wavelets will always increase the cardinality of a tree by $2$, while increasing the polynomial degree on a leaf $\lambda$ will instead increase the cardinality depending on the size of the set $\Upsilon(\lambda)$. Since this set has the form $\Upsilon(\lambda) = \{\mu \in \mathcal{T} : \lambda \succeq \mu \succeq \mu_\lambda\}$ for some $\mu_\lambda \preceq \lambda$ we can estimate
	\[1 \le \vert \Upsilon(\lambda) \vert \le \vert \{\mu \in \mathcal{T} : \lambda \succeq \mu \succeq \mathcal{R}\}\vert = \vert \lambda \vert +1.\]
	Refining $N_p = N$ times in polynomial degree on the node $\lambda=\mathcal{R}$ with $\vert \Upsilon(\mathcal{R}) \vert = 1$ directly gives the lower estimate in \eqref{T_N_cardinality}.
	
	For the upper bound we have to investigate how we can create the tree $T$ that maximizes  $\# T$ after $N$ refinement steps. In a single step, the largest increase in cardinality that is possible for a tree $(\mathcal{T},P_{\max})$ of depth $J = \max_{\lambda \in \mathcal{T}} |\lambda|$ by means of polynomial enrichment can occur if there exits a leaf $\lambda \in \mathcal{V}(\mathcal{T})$  with $\vert \lambda \vert = J$ and $\Upsilon(\lambda) 
	= \{\mu \in \mathcal{T} : \lambda \succeq \mu \succeq \mathcal{R}\}$. In this case polynomial enrichment of $\lambda$ will increase the cardinality of the quarklet tree by $\vert \Upsilon(\lambda) \vert = J +1$. On the other hand we avoid having many leaves on a high level since space refinement increases the cardinality only by $2$. Consequently the largest possible tree after $N$ refinement steps consists only of leaves and a single path to a leaf $\lambda$ on a high level with $\Upsilon(\lambda) 
	= \{\mu \in \mathcal{T} : \lambda \succeq \mu \succeq \mathcal{R}\}$ and polynomial enrichment is applied only on this leaf. This means we first have to employ $N_h$ steps of space refinement along this path such that we have $|\lambda| = N_h$. Then we refine the polynomial degree $N_p$-times on the leaf $\lambda$. The cardinality of such a tree is given by
	\begin{equation}\label{upper_bound_card}
	\#T = 2N_h +\left(N_h+1\right)\left(N-N_h \right) + 1,
	\end{equation}
	which has its maximum over $[0,N]$ in $N_h = \frac{N+1}{2}$. Inserting this into \eqref{upper_bound_card} yields the upper bound in \eqref{T_N_cardinality}.
	
	To prove \eqref{T_N_cardinality2} we start with the special case that $\max_{\lambda \in \mathcal{T}_N} |\lambda| = 0$. In this case the root $\mathcal{R}$ has been polynomially enriched $N$-times and we have $ \#T_N=N+1$. Now let $\max_{\lambda \in \mathcal{T}_N} |\lambda| \ge 1$ and let $\eta$ be a node of $\mathcal{T}_N$ with $|\eta| = \max_{\lambda \in \mathcal{T}_N} |\lambda|$. This especially implies $\eta \in \mathcal{V}(\mathcal{T}_N)$. Increasing the polynomial degree on the leaf $\eta$ raises the cardinality of the tree up to $|\eta|+1$. The total number of refinements is bounded by $N$ and hence the claim follows. 
\end{proof}

\subsection{The Trees $T_N$ are Near-Best}\label{subsec_alg_nearbest1}
Now we will prove that the trees produced by the algorithm {\bf NEARBEST\textunderscore TREE} are indeed near-best in the sense of \eqref{eq:nearbest}. We will split up this task into two substeps. At first we derive a lower bound for the best approximation error $\sigma_n$ in terms of the threshold parameter $q_N = q(\mathcal{R})$.
\begin{Lemma}\label{lemma:lower1}
Let $  n,N \in \mathbb{N}  $ with $n \le N$. Let $T^\star = (\mathcal{T}^\star, P_{\max}^\star)$ be the optimal quarklet tree of cardinality $n $ such that $\sigma_n = \mathcal{E}(T^\star)$. Let the tree $T_N = (\mathcal{T}_N, \mathcal{T}'_N) $ be produced by the algorithm {\bf NEARBEST\textunderscore TREE} and a subsequent trimming. Let us define the threshold parameter $q_N\coloneqq q(\mathcal{R})$ with respect to the tree $\mathcal{T}'_N$. Then it holds
	\begin{equation*}\label{eq:lower1}
		\sigma_n \ge q_N (N-n+1).
	\end{equation*}
\end{Lemma}

\begin{proof}
Let us consider the leaves $\lambda \in \mathcal{V}(\mathcal{T}^\star)$ and their orders $ P_{\max}^\star = \{p_{\textrm{max}}^\star(\lambda)\}_{\lambda \in \mathcal{V}(\mathcal{T}^\star)}$. In the case $r(\lambda, \mathcal{T}'_N) \le p_{\max}^\star(\lambda)$ we ignore the contribution of $e_{p_{\max}^\star(\lambda)}(\lambda)$ to the total error $\mathcal{E}(T^\star)$. That means we estimate
     \begin{equation}\label{eq:split}
     	\sigma_n = \mathcal{E}(T^\star) =  \sum_{\lambda \in \mathcal{V}(\mathcal{T}^\star)} e_{p_{\max}^\star(\lambda)}(\lambda) \ge \sum_{\substack{\lambda \in \mathcal{V}(\mathcal{T}^\star), \\ r(\lambda, \mathcal{T}'_N) > p_{\max}^\star(\lambda)}} e_{p_{\max}^\star(\lambda)}(\lambda).
     \end{equation} 
Now let $ k \in \mathbb{N}_0  $ with $ k \leq N  $. Then for the remaining terms of the sum which fulfill $r(\lambda, \mathcal{T}'_N) > p_{\max}^\star(\lambda)$ we consider the quantity $q_k \coloneqq q(\mathcal{R})$ at the stage $\mathcal{T}'_k$ of growing the tree $\mathcal{T}'_N$ at the last increase of $r(\lambda)$. In other words $\mathcal{T}'_k$ is the last tree in the sequence of trees $\mathcal{T}'_1,\ldots,\mathcal{T}'_N$ produced by the algorithm where a descendant of $\lambda$ is added. Notice that by \eqref{eq:2.3} and \eqref{eq:3.2} the quantities involved in the definition of $q(\lambda)$ are nonincreasing in the process of growing the tree. Therefore the numbers $q_k$ are decreasing with $k$ and we have $q_k \ge q_N$. By definition of $q$ and $s$ it follows that at this stage we have $s(\mathcal{R}) = s(\lambda) = s(\mu)$ for all $ \mu \preceq \lambda $. Now let $\mu$ be the parent of $\lambda$. Then we use \eqref{def_func_q1} to observe
\begin{align*}
	q(\mu) = \min \{q(\lambda), \tilde{E}_{r(\mu, \mathcal{T}'_k)}(\mu)\} \le q(\lambda).
\end{align*}
	By iterating this argument we obtain $q(\lambda) \ge q_k =q(\mathcal{R})$ and $\tilde{E}_j(\lambda) \ge q(\lambda) \ge q_N$ for $j=r(\lambda,\mathcal{T}'_N)-1$. Next we use the calculations provided in formula \eqref{eq:3.3} to find	
\begin{align*}
	\frac{1}{\tilde{E}_j(\lambda)} = \sum_{i = p_{\max}^\star(\lambda)+1}^j \frac{1}{E_i(\lambda)} + \frac{1}{\tilde{E}_{p_{\max}^\star( \lambda)}(\lambda)}.
\end{align*}
Then, by following the lines of \cite{Bin18}, see page 3355, we obtain the estimate
\begin{align*}
	E_{p_{\max}^\star(\lambda)}(\lambda) & \ge q_N(j-p_{\max}^\star(\lambda)+1).
\end{align*}
In a next step we use the definition of the errors $E_j$, see the formulas  \eqref{eq:3.1} and \eqref{eq:3.1a}, and recall $j=r(\lambda,\mathcal{T}'_N)-1$. Then we also get
	\begin{equation} \label{eq:3.6}
	e_{p_{\max}^\star(\lambda)}(\lambda) \ge E_{p_{\max}^\star( \lambda)}(\lambda) \ge q_N \max\{r(\lambda, \mathcal{T}'_N) - p_{\max}^\star( \lambda),0\}.
	\end{equation}
It remains to estimate the differences $r(\lambda, \mathcal{T}'_N) - p_{\max}^\star( \lambda)$ for all leaves $\lambda \in \mathcal{V}(\mathcal{T}^\star)$. With $p_{\max}^\star(\lambda) $ for $\lambda \in \mathcal{T}^\star \backslash \mathcal{V}(\mathcal{T}^\star)$ being induced by $P_{\max}$ we observe
	\begin{equation}\label{eq:A4}
	\sum_{\lambda \in \mathcal{V}(\mathcal{T}^\star)} \max\{r(\lambda, \mathcal{T}'_N)- p_{\max}^\star(\lambda),0\}\ge\sum_{\lambda \in \mathcal{V}(\mathcal{T}^\star \cap \mathcal{T}'_N)} r(\lambda, \mathcal{T}'_N)- p_{\max}^\star(\lambda). 
	\end{equation}
To further estimate the right side of \eqref{eq:A4} we need two intermediate estimates. On the one hand we have
		\begin{align*}
	\sum_{\lambda \in \mathcal{V}(\mathcal{T}^\star \cap \mathcal{T}'_N)} p_{\max}^\star(\lambda) \le \sum_{\lambda \in \mathcal{V}(\mathcal{T}^\star)} p_{\max}^\star(\lambda) \le  n - |\mathcal{T}^\star|,
	\end{align*}
	while on the other hand it holds
	\begin{align*} 
	    \sum_{\lambda \in \mathcal{V}(\mathcal{T}^\star \cap \mathcal{T}'_N)} r(\lambda, \mathcal{T}'_N) = N - \frac{|\mathcal{T}^\star \cap \mathcal{T}'_N|-1}{2} \ge N - |\mathcal{T}^\star|+1.
	\end{align*}
	We combine these estimates with \eqref{eq:A4} to obtain
	\begin{equation}\label{eq:3.7}
		\sum_{\lambda \in \mathcal{V}(\mathcal{T}^\star)} \max\{r(\lambda, \mathcal{T}'_N)- p_{\max}^\star(\lambda),0\} \ge N-n+1.
	\end{equation}
Finally, a combination of \eqref{eq:split}, \eqref{eq:3.6} and \eqref{eq:3.7} yields
	\begin{align}
	\sigma_n &= \sum_{\lambda \in \mathcal{V}(\mathcal{T}^\star)} e_{p_{\max}^\star(\lambda)}(\lambda) \ge q_N \sum_{\lambda \in \mathcal{V}(\mathcal{T}^\star)}\max\{r(\lambda, \mathcal{T}'_N)- p_{\max}^\star(\lambda),0 \}\notag \\
	&\ge q_N(N-n+1). \notag
	\end{align}
\end{proof}

In a next step we want to provide an upper bound for the global error with respect to the threshold parameter $q_N$.

\begin{Lemma}\label{lemma:upper}
Let $  N \in \mathbb{N}  $ and let the tree $T_N = (\mathcal{T}_N, \mathcal{T}'_N) $ be produced by the algorithm {\bf NEARBEST\textunderscore TREE} and a subsequent trimming. Let $q_N = q(\mathcal{R})$ with respect to the tree $\mathcal{T}'_N$. Then for the global error there is the upper bound
	\begin{equation}\label{lem_upper:estimate}
		\mathcal{E}(T_N) \le  q_N \left(2N+1\right).
	\end{equation}
\end{Lemma}

\begin{proof}
In general this result can be proved with similar methods as described in \cite{Bin18}, see pages 3355-3356. Nevertheless, some modifications are necessary. We consider the tree $ \mathcal{T}'_N  $. Let $L$ denote the set of nodes $\lambda \in \mathcal{T}'_N  $ for which we have $q(\lambda) = \tilde{E}_{r(\lambda, \mathcal{T}'_N)}(\lambda)$. Moreover let $Q$ be the maximal subtree of $\mathcal{T}'_N$ emanating from $\mathcal{R}$ such that it holds $L \cap Q = \mathcal{V}(Q)$. Consequently there is no inner node of $Q$ in the set $L$. As in \cite{Bin18} we observe that $Q$ is a complete tree.  Moreover for $ \lambda \in \mathcal{V}(Q)  $ we get
	\begin{equation*} \label{eq:B3}
		\tilde{E}_{r(\lambda, \mathcal{T}'_N)}(\lambda) = q(\lambda) \le  q(\mathcal{R}) = q_N.
	\end{equation*}
Now let $\lambda \in \mathcal{V}(Q)$ and put $ r( \lambda, \mathcal{T}'_N) = j$. Then with similar arguments as in \cite{Bin18}, see especially formula (3.9) on page 3356, we find 
\begin{equation} \label{eq:A}
	E_{j}(\lambda) = \tilde{E}_{j}(\lambda)\left( \sum_{k=0}^{j} \frac{E_{j}(\lambda)}{E_k(\lambda)}+\sum_{\mu \prec \lambda} \frac{E_j(\lambda)}{e(\mu)}\right) \le q_N \left( j+1 + \sum_{\mu \prec \lambda}\frac{e(\lambda)}{e(\mu)} \right).
\end{equation}
Notice that for a leaf $\lambda \in \mathcal{V}(Q)$ two different cases can show up. On the one hand it is possible that there exists a $\nu \in \mathcal{V}(\mathcal{T}_N)$ with $\lambda \preceq \nu$. Then we can write
	\begin{equation}\label{eq:desc}
		E_{r(\lambda, \mathcal{T}'_N)}(\lambda) = \sum_{\eta \in \mathcal{V}(\mathcal{J}_{\lambda} \cap \mathcal{T}_N)} E_{r(\eta, \mathcal{T}'_N)}(\eta) = \sum_{\eta \in \mathcal{V}(\mathcal{J}_{\lambda} \cap \mathcal{T}_N)} e_{r(\eta, \mathcal{T}'_N)}(\eta).
	\end{equation}
On the other hand it might be possible that there is a $\nu \in \mathcal{V}(\mathcal{T}_N)$ such that $\lambda \succ \nu$. In that case we observe
	\begin{equation}\label{eq:prec}
		 e_{r(\nu, \mathcal{T}'_N)}(\nu) = E_{r(\nu, \mathcal{T}'_N)}(\nu) \le  \sum_{\eta \in \mathcal{V}(\mathcal{T}'_{\nu} \cap Q)} E_{r(\eta, \mathcal{T}'_N)}(\eta). 
	\end{equation}
Consequently we can split the leaves $\mathcal{V}(\mathcal{T}_N)$ in two sets corresponding to the two cases. Then a combination of \eqref{eq:desc} and  \eqref{eq:prec} yields 
	\begin{equation*}\label{eq:B}
		\mathcal{E}(T_N) = \sum_{\lambda \in \mathcal{V}(\mathcal{T}_N)} e_{r(\lambda, \mathcal{T}'_N)}(\lambda) \le \sum_{\lambda \in \mathcal{V}(Q)} E_{r(\lambda, \mathcal{T}'_N)}(\lambda).
	\end{equation*}
Next we plug in \eqref{eq:A}. Then we find
	\begin{align}
	\mathcal{E}(T_N) &\le q_N \left(\sum_{\lambda \in \mathcal{V}(Q)} \left( r(\lambda, \mathcal{T}'_N) + 1\right)  + \sum_{\lambda \in \mathcal{V}(Q)} \sum_{\mu \prec \lambda} \frac{e(\lambda)}{e(\mu)}  \right) \notag\\
	& = q_N \left( \sum_{\lambda \in \mathcal{V}(Q)} (r(\lambda, \mathcal{T}'_N) +1)  +  \sum_{\mu\in (Q \backslash \mathcal{V}(Q))} \frac{\sum\limits_{ \lambda \in \mathcal{V}(\mathcal{J}_{\mu} \cap Q)} e(\lambda)}{e(\mu)}  \right). \label{eq:global_err_est}
	\end{align}
Here in the last step we changed the order of summation. Notice that $(\ref{eq:1.2})$  implies the weak subadditivity property 
\begin{equation*} \label{eq:weakSub}
e(\lambda) \ge  \sum_{\eta \in \mathcal{V}(\mathcal{J}_{\lambda} \cap \mathcal{T})} e(\eta)
\end{equation*} 
for all trees $\mathcal{T}$. This especially implies that each fraction in the second sum in \eqref{eq:global_err_est} is bounded from above by $1$. Consequently we get
	\begin{equation}\label{eq:A2}
		\mathcal{E}(T_N)  \le q_N  \left(|\mathcal{V}(Q)|+\sum_{\lambda \in \mathcal{V}(Q)} r(\lambda, \mathcal{T}'_N)   +  |Q\backslash \mathcal{V}(Q)|\right) .
	\end{equation}
Observe that the number of nodes in $\mathcal{T}'_N$ can be rewritten as
	\begin{equation*}\label{eq:B2}
		|\mathcal{T}'_N| = |Q| + 2\sum_{\lambda \in \mathcal{V}(Q)} r(\lambda, \mathcal{T}'_N).
	\end{equation*}
We insert this into \eqref{eq:A2} to obtain
	\begin{align*}
		\mathcal{E}(T_N) & \le q_N \left( |Q|+ \frac{|\mathcal{T}'_N|-|Q|}{2} \right).
	\end{align*}
Let us recall that $ |Q| \le |\mathcal{T}'_N| = 2N+1$. Then we finally get 
	\begin{align*}
	 \mathcal{E}(T_N)  \le q_N \left( 2N+1\right)  .
	 \end{align*}
This is what we stated in \eqref{lem_upper:estimate}.
\end{proof}

Now we are well-prepared to prove that the algorithm  {\bf NEARBEST\textunderscore TREE} provides a quarklet tree which is a near-best approximation in the sense of \eqref{eq:nearbest}.

\begin{Theorem}\label{theorem:1}
Let $ n,N \in \mathbb{N} $ with $n\le N$ and let $ e_{p}(\lambda)  $ be local errors that fulfill \eqref{eq:1.2} and \eqref{eq:1.3}. Then {\bf NEARBEST\textunderscore TREE} with a subsequent trimming finds a quarklet tree $T_N = (\mathcal{T}_N, \mathcal{T}'_N)$ such that the corresponding quarklet approximation is near-best in the sense
	\begin{equation}\label{eq:result1}
		\mathcal{E}(T_N) \le \frac{2N+1}{N-n+1}\sigma_n.
	\end{equation}
\end{Theorem}

\begin{proof}
Formula \eqref{eq:result1} follows by combining Lemma \ref{lemma:lower1} and Lemma \ref{lemma:upper}.
\end{proof}

\begin{Remark} \label{remark:interpretation}
	Let $M \ge 2$ be a natural number and let us ignore the fact that our algorithm can only perform integer steps. Then the first part of Lemma  \ref{lemma:T_N_cardinality} implies that we can run at most $N_1 = 2\sqrt{M+1}-3$ steps of our algorithm while guaranteeing that the resulting tree $T_{N_1}$ fulfills $\#T_{N_1} \le M$. Now let $n_1 = \frac{N_1}{2} $. Using the trivial estimate
\begin{align*}
   \frac{2N_1+1}{N_1-n_1+1} = \frac{2\left(2\sqrt{M+1}-3\right)+1}{2\sqrt{M+1}-3- \sqrt{M+1}+ \frac52} \le 4,
\end{align*}
formula \eqref{eq:result1} then becomes
\begin{align*}
	\mathcal{E}(T_{N_1}) \le 4\sigma_{ \sqrt{M+1}-\frac32}.
\end{align*}
In particular, this implies
\[	\mathcal{E}(T_{N_1}) \le C_1\sigma_{c_1M^{1/2}},\]
with independent constants $C_1>0$ and $c_1 \in (0,1]$. In practice, the global maximal refinement level of the quarklets is usually bounded by some $j_{\textrm{MAX}} \in \mathbb{N}$. Then \eqref{T_N_cardinality2} implies $\# T_N \le C N$, where the constant $C>0$ depends on the maximal refinement level $j_{\textrm{MAX}} $. Therefore we can run $N_2 =  C^{-1}M$ steps of our algorithm while guaranteeing that the resulting tree $T_{N_2}$ fulfills $\# T_{N_2}\le M$. With $n_2 = \frac{N_2}{2}$ the trivial estimate
	\begin{align*}
	\frac{2N_2+1}{N_2-n_2+1} = \frac{2 C^{-1}M+1}{ C^{-1}M- C^{-1}\frac{M}{2}+1} \le 4
	\end{align*}
	then yields that \eqref{eq:result1} becomes
	\begin{align*}
	\mathcal{E}(T_{N_2}) \le 4\sigma_{\frac{M}{2C}}.
	\end{align*}
	In particular, this implies
	\[	\mathcal{E}(T_{N_2}) \le C_2\sigma_{c_2M},\]
	with an independent constant $C_2>0$ and a constant $c_2 \in (0,1]$ that depends on the maximal refinement level.
\end{Remark}

\section{Practical Realization} \label{sec_practice}

\subsection{Approximation in $\ell_2$}\label{sec:coeff_approach}

In this section we describe one possible way how to apply {\bf NEARBEST\textunderscore TREE} to functions $ f \in L_2(I)  $. To avoid additional notation and technical difficulties when dealing with boundary adapted quarks and quarklets we restrict us to the Haar quarklet case $m = \tilde{m} = 1$ where we have $\supp \psi_p = \supp \varphi_p = [0,1]$ and $\supp \psi_{p,j,k} = 2^{-j}[k,k+1]$. The infinite set of the quarklet indices corresponding to quarklets which intersect with $I$ is consequently given by 
\[ \Lambda_I = \{(p,j,k_j):p\in \mathbb{N}_0, j \in \mathbb{N}_0\cup\{-1\}, 0 \le k_j \le k_{j,\max}\coloneqq \max\{0,2^j-1\} \}.\]
For the corresponding set of wavelet indices we write
\[ \Lambda_{I,0} = \{(j,k_j): j \in \mathbb{N}_0\cup\{-1\}, 0 \le k_j \le k_{j,\max} \}.\]
This fits into the setting of dyadic intervals, see Section \ref{sec:TreeIndices}, with one small exception. We set $\mathcal{R} = (0,0)$ and also assign the indices with $j=-1$ to the root. This is sensible because we have no direct parent child relation between the quarklets $\psi_{p,0,k}$ and the quarks $\psi_{p,-1,k} = \varphi_p(\cdot-k)$ since $\supp \psi_{p,-1,k} = \supp \psi_{p,0,k}$ and $k_{-1,\max} = k_{0,\max} = 0$. This means if we consider a tree of quarklet indices $T = (\mathcal{T},P_{\max})$ and increase the maximal polynomial degree  $p = p_{\max}(\mathcal{R})$ on the root by one, we instead add two indices to the tree, namely the quark index $(p+1,-1,0)$ and the quarklet index $(p+1,0,0)$. Furthermore, when we sum over all ancestors of an index $\lambda \in \mathcal{T}$ we also include the quark index.

\begin{Remark}
	When working with quarklets of order $m \ge  2$ we have  $\supp \psi_p \nsubseteq [0,1]$. In this case one usually starts with a minimal level $j_{\min}>0$ and employs additional boundary quarklets, see \cite{bib:DFK18} for details. The infinite index set in this situation consists of multiple nodes on the level $j_{\min}$. This means we have $R>1$ roots. To deal with this situation one can employ the strategy stated in \cite{bib:CNSV14}. In case $R$ is a power of $2$, we unify the existing roots by forming pairs of them and creating a new parent for each pair. This procedure is then repeated until only one root remains. Otherwise we create up to $\lceil \log_2(R) \rceil -1$ empty nodes before the unification process. For more details see Remark 6 in \cite{bib:CNSV14}.
\end{Remark}

Now we turn our focus to the precise definition of the local errors $e_p(\lambda)$. In the previous sections we made some direct and implicit assumptions which we recall here for clarity:
\begin{itemize}
	\item There is subadditivity of the error for the lowest order, see \eqref{eq:1.2}.
	\item There is a reduction of the error when the maximal polynomial degree increases, see \eqref{eq:1.3}.
	\item The global error $\mathcal{E}(T)$ is equivalent to the current approximation error, see Remark \ref{rem_glob_err}.
\end{itemize}
It will be our aim to derive a suitable definition of these local errors, which we then can apply in practice in Section \ref{numerical_results}. We start by recalling that every $ f \in L_2(I)  $ can (nonuniquely) be written as  
\begin{equation}\label{rep_L2_gen}
	f =  \sum_{(p,j,k) \in  \Lambda_I}  c_{p,j,k} w_p \psi_{p,j,k}
\end{equation}
with a coefficient sequence $ \mathbf{c} = \{ c_{p,j,k}  \}_{(p,j,k) \in  \Lambda_I} $. We are now ready to provide a suitable definition of the local errors. 

\begin{Definition}\label{def_locerr_L2_new}
	Let $ f \in L_2(I)  $ be given in the form \eqref{rep_L2_gen}. Then for each node $  \lambda \in \Lambda_{I,0}$ and $p \in \mathbb{N}_0$ we define the local errors $  e_p(\lambda) $ via 
	\begin{equation}\label{eq:e_p_coeffs_new}
	e_p(\lambda) \coloneqq \sum_{(i,\ell) \in \Upsilon(\lambda)} \sum_{q > p} \vert c_{q,i,\ell} \vert^2 + \sum_{(i,\ell)  \succ \lambda} \sum_{q \ge 0} \vert c_{q,i,\ell} \vert^2.
	\end{equation}
\end{Definition}

The first sum in \eqref{eq:e_p_coeffs_new} refers to $  \lambda $ and a subset of its ancestors, whereby the polynomial degree is greater than the maximal degree $p$. The second sum gathers all descendants of $ \lambda $ for all possible polynomial degrees. 
\begin{Remark}
	Definition \ref{def_locerr_L2_new} is inspired by tree approximations from adaptive wavelet schemes. In particular, in the context of nonlinear variational problems approximation with tree structured index sets based on the size of the wavelet expansion coefficients has already been successfully employed, see for example \cite{bib:CDD03} and \cite{kappei2011adaptive}.
\end{Remark}

In the following, we want to check that  the essential conditions \eqref{eq:1.2} and \eqref{eq:1.3} mentioned in Section \ref{subsec_localerr} are fulfilled.

\begin{Lemma}\label{lem_loc_err_prop1_new}
	Let $ f \in L_2(I)  $ be given in the form \eqref{rep_L2_gen}. Let $\lambda \in \Lambda_{I,0}$ and $p \in \mathbb{N}_0$.  Then the local errors $ e_p(\lambda)  $ from Definition \ref{def_locerr_L2_new} satisfy the properties \eqref{eq:1.2} and \eqref{eq:1.3}.
\end{Lemma}

\begin{proof}
	At first we prove that \eqref{eq:1.2} holds. Therefore let $p=0$ and $ \eta_1, \eta_2 $ be the children of $ \lambda  $. Notice that by definition of the sets $\Upsilon$ we have $ \Upsilon(\lambda) \cup \{\eta_1, \eta_2\} = \Upsilon(\eta_1) \cup \Upsilon(\eta_2) $, see \eqref{eq:Upsilon} and below. Then for the local error of the lowest order concerning $ \lambda  $ we find
	\begin{align*}
	e_0(\lambda) &= \sum_{(i,\ell) \in \Upsilon(\lambda)} \sum_{q > 0} \vert c_{q,i,\ell} \vert^2 + \sum_{(i,\ell) \succ \lambda} \sum_{q \ge 0} \vert c_{q,i,\ell} \vert^2 \\
	&\ge \sum_{(i,\ell) \in \Upsilon(\lambda)} \sum_{q > 0} \vert c_{q,i,\ell} \vert^2 + \sum_{(i,\ell) \in \{\eta_1, \eta_2 \} } \sum_{q > 0} \vert c_{q,i,\ell} \vert^2 + \sum_{(i,\ell) \succ\eta \in \{\eta_1, \eta_2\}} \sum_{q \ge 0} \vert c_{q,i,\ell} \vert^2\\
	&= \sum_{(i,\ell) \in \Upsilon(\eta_1)} \sum_{q > 0} \vert c_{q,i,\ell} \vert^2 + \sum_{(i,\ell) \in \Upsilon(\eta_2)} \sum_{q > 0} \vert c_{q,i,\ell} \vert^2 + \sum_{(i,\ell) \succ\eta \in \{\eta_1, \eta_2\}} \sum_{q \ge 0} \vert c_{q,i,\ell} \vert^2\\
	& = e_0(\eta_1) +  e_0(\eta_2).
	\end{align*}
	Property \eqref{eq:1.3} follows directly from Definition \ref{def_locerr_L2_new}. 
\end{proof}

Now we are prepared to run the algorithm {\bf NEARBEST\textunderscore TREE} with the subsequent trimming routine {\bf TRIM}. We will show that the global error indeed describes the quality of the approximation.

\begin{Lemma}\label{lem_L2_globerr}
	Let $ f \in L_2(I)  $ be given in the form \eqref{rep_L2_gen}. Let the local errors $ e_p(\lambda)   $ be defined as in Definition \ref{def_locerr_L2_new}.  For $ N \in \mathbb{N} $ by $T_{N} = (\mathcal{T}_N, P_{\max})$ we denote the quarklet tree resulting out of the algorithm {\bf NEARBEST\textunderscore TREE} with a subsequent trimming. The corresponding approximation $ f_{T_N} $ is given by
	\begin{equation}\label{approx_fN_L2}
	f_{T_N} = \sum_{(p,j,k) \in T_{N} \subset \Lambda_I }  c_{p,j,k} w_p\psi_{p,j,k} .
	\end{equation} 
	Then there exists a constant $ C > 0 $ independent of $ f $ and $ N $ such that for the global error we observe
	\begin{align}\label{eq:reliable}
	 \Vert f - f_{T_N} \Vert_{L_2(I)}^2 \leq C \mathcal{E}(T_{N}) .
	\end{align}
\end{Lemma} 

\begin{proof}
	To prove this result in a first step we use the definition of the global error, see \eqref{eq_glob_err}, in combination with Definition \ref{def_locerr_L2_new}. Then we find 
	\begin{align*}
	\mathcal{E}(T_{N})  = \sum_{\lambda \in \mathcal{V}(\mathcal{T}_N)} e_{p_{\max}(\lambda)}(\lambda)  = \sum_{\lambda \in \mathcal{V}(\mathcal{T}_{N})} \left( \sum_{(i,\ell)  \in \Upsilon(\lambda)} \sum_{q > p_{\max}(\lambda)} \vert c_{q,i,\ell} \vert^2 + \sum_{(i,\ell)  \succ \lambda} \sum_{q \ge 0} \vert c_{q,i,\ell} \vert^2 \right)   .
	\end{align*} 
	Recall that by definition of the sets $\Upsilon(\lambda)$ we have $ \bigcup_{ \lambda \in \mathcal{V}(\mathcal{T}_N)} \Upsilon(\lambda) = \mathcal{T}_N$. Using this observation we can also write 
	\begin{align*}
	\mathcal{E}(T_{N}) =  \sum_{\lambda = (j,k) \in \mathcal{T}_{N}} \sum_{ q > p_{\max}(\lambda)  }  \vert c_{q,j,k} \vert^2
	+ \sum_{(i,\ell) \succ \lambda \in \mathcal{V}(\mathcal{T}_{N})} \sum_{q \geq 0 }  \vert c_{q,i,\ell} \vert^2  .
	\end{align*} 
    Let us have a closer look at this expression. The first sum collects all quarklet coefficients whose corresponding wavelet indices are nodes of the tree $ \mathcal{T}_{N} $. But nevertheless these coefficients do not belong to the tree $ T_{N}  $ due to the polynomial degree of the corresponding quarklets. The second sum collects all coefficients that do not belong to  $ T_{N}  $ since their corresponding indices are descendants of the leaves of $ \mathcal{T}_{N} $. Therefore we can write
	\begin{align*}
	\mathcal{E}(T_{N}) =   \sum_{(p,j,k) \in \Lambda_I}  \vert c_{p,j,k} \vert^2 -    \sum_{(p,j,k) \in T_{N} \subset \Lambda_I}  \vert c_{p,j,k} \vert^2.
	\end{align*}
	Next we use the frame property from Theorem \ref{L_2_frame} together with the lower estimate in \eqref{eq:frame_charac} to obtain
	\begin{align*}
	\Big \Vert   \sum_{(p,j,k) \in \Lambda_I}  c_{p,j,k} w_p \psi_{p,j,k} -    \sum_{(p,j,k) \in T_{N} \subset \Lambda_I}  c_{p,j,k} w_p \psi_{p,j,k} \Big \Vert_{L_2(I)}^2 \lesssim \mathcal{E}(T_{N}).
	\end{align*} 
	Recall that $ f \in L_2(I) $ has the form \eqref{rep_L2_gen}. Then with \eqref{approx_fN_L2} we finally get 
	\begin{align*}
	\Vert   f - f_{T_N}  \Vert_{L_2(I)}^2 \lesssim \mathcal{E}(T_{N}) .
	\end{align*}
	This completes the proof.
\end{proof}  

Let us remark that an error estimator fulfilling \eqref{eq:reliable} is called \emph{reliable}. Now we are well prepared to describe the quality of the approximation $  f_{T_N}    $ concerning $ f $. For that purpose we work with the best tree approximation error $\sigma_n$, see Definition \ref{def_bestappr_err}, and apply Theorem \ref{theorem:1}. Then we obtain the following result. 

\begin{Lemma}\label{lem_L2_bestappr}
	Let $ f \in L_2(I)  $ be given in the form \eqref{rep_L2_gen} and let the local errors $ e_p(\lambda)   $ be defined as in Definition \ref{def_locerr_L2_new}.  For $ N \in \mathbb{N} $ by $T_N = (\mathcal{T}_{N},P_{\max})$ we denote the quarklet tree resulting out of the algorithm {\bf NEARBEST\textunderscore TREE} with a subsequent trimming and cardinality $\#T_{N} \lesssim N $ whereby the constant depends on the maximal level of the quarklet indices in the tree. The corresponding approximation $ f_{T_N} $ is given via
	\begin{equation*}\label{approx_fN_L2.2}
	f_{T_N} = \sum_{(p,j,k) \in T_{N}}  c_{p,j,k} w_p \psi_{p,j,k}.
	\end{equation*} 
	Let $ n \le N  $. Then there exists a constant $ C > 0 $ independent of $ f $, $ N $ and $ n $ such that 
	\begin{align*}
	\Vert f - f_{T_N} \Vert_{L_2(I)}^2 \leq C \frac{2N+1}{N-n+1}\sigma_n  .
	\end{align*}
\end{Lemma} 

\begin{proof}
	To prove this result at first we use Lemma \ref{lem_L2_globerr}. Notice that all conditions from there are fulfilled. We obtain 
	\begin{align*}
	\Vert f - f_{T_N}  \Vert_{L_2(I)}^2 \leq C \mathcal{E}(T_{N}) .
	\end{align*} 
	Now we want to apply Theorem \ref{theorem:1}. For that purpose recall that the local errors $ e_p(\lambda)  $ from Definition \ref{def_locerr_L2_new} satisfy the properties \eqref{eq:1.2} and \eqref{eq:1.3}, see Lemma \ref{lem_loc_err_prop1_new}. Consequently we get 
	\begin{align*}
	\Vert f - f_{T_N} \Vert_{L_2(I)}^2 \leq C \frac{2N+1}{N-n+1}\sigma_n  .
	\end{align*}
	This completes the proof. 
\end{proof}

\begin{Remark}\label{rem_fs_general}
	The results from this section can be transferred from the $L_2$-setting to more sophisticated function spaces like Sobolev spaces $H^s$, Besov spaces $B_{r,q}^s$ or Triebel-Lizorkin spaces $F_{r,q}^s$. As we have seen in the proof of Lemma $\ref{lem_L2_globerr}$ the important aspect is the equivalence between the norm of the function space and the sequence norm of the quarklet coefficients. We refer to \cite{bib:DKR17},\cite{bib:sieber2020adaptive} and \cite{HoDa2021} for the frame property of the quarklet system in different function spaces. However in the case of Besov spaces and Triebel-Lizorkin spaces the proofs become much more technical. Moreover in those cases some additional conditions concerning the parameters $s,r$ and $q$ will become necessary.
\end{Remark}

\begin{Remark}\label{rem:coeff_calc}
	So far we have not commented on the task of actually finding a representation of $f \in L_2(I)$ as in \eqref{rep_L2_gen}. In some applications our function is already given as a linear combination of quarklets. For example this is the case if $f$ stems from prior computations, so $f$ could be the current approximation of the unknown solution of a linear elliptic operator equation discretized with a quarklet frame. However there are also situations where this is not the case. Then we have to look for a suitable coefficient sequence $\{c_\lambda\}_{\lambda \in \Lambda_I}$. As stated earlier, this representation is not necessarily unique. Therefore the results from applying the machinery presented in this section will depend on the choice of the coefficient sequence. In practice, we always work with a finite subset $ \overline{\Lambda}_I \subset \Lambda_I$. For example, we can truncate $\Lambda_I$ at a uniform maximal refinement level and polynomial degree.  Then one way to calculate the coefficient sequence $\textbf{c} \in \ell_2(\overline{\Lambda}_I)$ is given by solving the matrix-vector equation
	\begin{equation}\label{eq:Gc=F}
	\mathbf{Gc}=\mathbf{b},
	\end{equation}
	where $\mathbf{b} \coloneqq \left( \langle f, w_\lambda \psi_\lambda \rangle_{L_2(I)}\right)_{\lambda \in \overline{\Lambda}_I} $ and $\mathbf{G} \coloneqq \left( \langle w_\lambda\psi_\lambda, w_\mu\psi_\mu \rangle_{L_2(I)}\right)_{\lambda,\mu \in \overline{\Lambda}_I} $ denotes the Gramian. However, the Gramian matrix has a non-trivial kernel due to the redundancy in the quarklet system. Therefore \eqref{eq:Gc=F} is not uniquely solvable. Nonetheless, classical iterative schemes like the damped Richardson iteration
	\[\mathbf{c}^{(j+1)} \coloneqq \mathbf{c}^{(j)}+\omega(\mathbf{b}-\mathbf{G}\mathbf{c}^{(j)}), \quad 0 < \omega < \frac{2}{||\mathbf{G}||_2}, \quad j = 0,1,\ldots\]
	or variations thereof can still be applied in a numerically stable way. Unfortunately, coefficients derived this way sometimes have a poor qualitative behavior in the sense that there are coefficients on a high level with a large modulus in regions where $f$ is considered to be `smooth'. This is counter intuitive since smooth parts should be well resolved by quarklets of a low level and with high polynomial degrees. If we now apply {\bf NEARBEST\textunderscore TREE} with the local errors as in \eqref{eq:e_p_coeffs_new} this choice of coefficients leads to many refinements towards these large coefficients. This in turn negatively impacts the resulting convergence rate. From the practical point of view we are also more interested in sparsity of the representation than finding the solution with the smallest norm. For these reasons we proceed iteratively. We start with a small index set $\tilde{\Lambda} \subset \overline{\Lambda}_I$ and solve \eqref{eq:Gc=F} with respect to $\tilde{\Lambda}$ instead of $\overline{\Lambda}_I$. Then we apply the zero extension operator on the solution to end up with a sequence $\mathbf{c} \in \ell(\overline{\Lambda}_I)$. Based on the information provided from the residual $\mathbf{Gc-b}$ we then add additional wavelet or quarklet indices to $\tilde{\Lambda}$. This procedure is repeated until the residual error $\Vert \mathbf{Gc-b} \Vert_2 $ is smaller than a given tolerance. By proceeding this way we end up with coefficient sequences with a better qualitative behavior at the expense of solving the matrix-vector equation \eqref{eq:Gc=F} only up to a certain accuracy. We still have to solve matrix-vector equations where the Gramian matrix has a non-trivial kernel. However these systems are much smaller than the whole truncated index set and can be handled easier by either classical iterative schemes or specialized software packages. This procedure is not completely satisfactory from the theoretical point of view. Finding a more reliable method to calculate suitable quarklet coefficients, e.g. in terms of an appropriate dual frame is the subject of further research.
\end{Remark}
	
\subsection{Numerical Experiments}\label{numerical_results}

In this section we test the algorithm {\bf NEARBEST\textunderscore TREE} with the local errors from Section \ref{sec:coeff_approach}. In the long run, we intend to use our algorithm as a building block in the design of  adaptive numerical methods to solve elliptic operator equations. Therefore the test examples are chosen as (models of) typical solutions of partial differential equations where we expect that adaptive schemes outperform classical uniform schemes.  In particular, for second order elliptic boundary value problems on polygonal domains with re-entrant corners it is well known that the exact solution consists of a regular and a singular part. A simple model for edge singularities is the univariate function  $x^\alpha$  with $\alpha > \frac12$, see e.g. \cite{bib:BS94}. It was shown in \cite{bib:DRS19} that this function can be directly approximated from the span of the quarklet system $\Psi$ at inverse-exponential rates. To be precise, the rate in \cite{bib:DRS19} for the approximation in $L_2$ is given by $e^{-2 \ln(2)n^{1/5}}$, whereby $n \in \mathbb{N}$ denotes the number of degrees of freedom. Therefore, in our setting exponential convergence refers to decay rates of
the form $e^{-\beta n^{\gamma}}$ for some $\beta, \gamma > 0$. This will serve as a benchmark for the approximations provided by the adaptive algorithm {\bf NEARBEST\textunderscore TREE}.

\begin{figure}[!ht] 
	\centering
	\begin{tikzpicture}[scale=0.9, 
	level 1/.style={sibling distance=70mm},
	level 2/.style={sibling distance=40mm},
	level 3/.style={sibling distance=20mm},
	level 4/.style={sibling distance=10mm}]
	\node [circle,draw, ultra thick,label=above:{$\mathcal{R} = (0,0)$}] (z){}
	child[draw,ultra thick] {node [circle,draw,label=left:{$(1,0)$}] (a) {}
		child[draw] {node [circle,draw,label=left:{$(2,0)$}] (b) {}
			child[draw]  {node [circle,draw, solid, label=left:{$(3,0)$}] (o){}}
			child[draw, dashed]  {node [circle,draw, solid, label=left:{$(3,1)$}] {} 
				child[draw, solid] {node [circle,draw, solid, label=left:{$(4,2)$}] (o) {}}
				child[draw, dashed] {node  [circle,draw, solid, label=right:{$(4,3)$}] (c){}}}}
		child[draw, dashed]  {node [circle,draw,solid, label=left:{$(2,1)$}] (g) {}
			child[draw, solid]  {node [circle,draw,solid, label=left:{$(3,2)$}] (o) {}}
			child[draw]  {node  [circle,draw,solid, label=left:{$(3,3)$}] (c){} }
		}
	}
	child[draw=black, dashed, ultra thick] {node [circle,draw, solid, label=right:{$(1,1)$}] (j) {}
		child[draw, solid] {node [circle,draw, solid, label=right:{$(2,2)$}] (k) {}
		}
		child[draw=black, dashed] {node [circle,draw, solid, label=right:{$(2,3)$}] (l) {}
			child[draw, solid] {node [circle,draw, solid, label=right:{$(3,6)$}] (o) {} 				child[draw, solid] {node [circle,draw, solid, label=left:{$(4,12)$}] (f) {}}
				child[draw, dashed] {node  [circle,draw, solid, label=right:{$(4,13)$}] (d){}}}
			child[draw, dashed] {node  [circle,draw, solid, label=right:{$(3,7)$}] (c){}  }
		}
	};
	\end{tikzpicture}
	\caption{Example of a tree $\mathcal{T}$ and the sets $\Upsilon$. Each node corresponds to a wavelet index $\lambda \in \mathcal{T}$. The solid lines represent the sets $\Upsilon(\lambda)$, $\lambda \in \mathcal{V}(\mathcal{T})$, with $\mu_\lambda$ set as in \eqref{eq:mu_lambda}.}
	\label{Fig:hp-grid}
\end{figure}
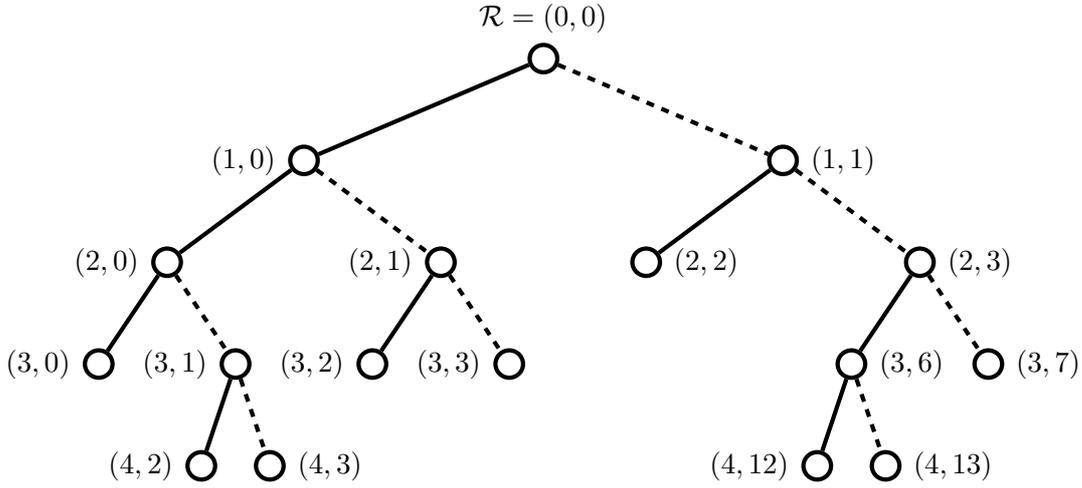

In all our experiments we consider the unit interval $I = [0,1]$ and the corresponding Haar quarklet index set $\Lambda_I$ from Section \ref{sec:coeff_approach} truncated at the uniform maximal refinement level $j_{\textrm{MAX}} = 10$ and maximal polynomial degree $p_{\textrm{MAX}} = 5$ denoted by $\overline{\Lambda}_I$. Furthermore the local errors $e_p(\lambda)$ are given by \eqref{eq:e_p_coeffs_new}.  That means we let $f \in L_2(I)$ be given in the form  \eqref{rep_L2_gen} with suitable coefficients $  \mathbf{c} = \{ c_{p,j,k}  \}_{(p,j,k) \in  \overline{\Lambda}_I} $, see Remark \ref{rem:coeff_calc}.
For the numerical experiments we also have to specify how the sets $\Upsilon(\lambda)$ look like. To this end we classify each wavelet index $(j,k) = \lambda \in \Lambda_{I,0}$ into one of two distinct groups. If $k \in 2\mathbb{N}_0$ we call $\lambda$ a \emph{left} node. If otherwise $k \in 2\mathbb{N}_0+1$ we call it a \emph{right} node. For a given node $\lambda \in \Lambda_{I,0}$ we then fix a specific ancestor $\mu_\lambda$ by setting
\begin{equation}\label{eq:mu_lambda}
	\mu_\lambda = \begin{cases}
	\mathcal{R}, & \textrm{all }  \mu \preceq \lambda  \textrm{ are left nodes}, \\
	\underset{\mu \preceq \lambda}{\arg \max} \{ \vert \mu \vert : \mu \textrm{ is a right node}\}, & \textrm{else}.
	\end{cases}
\end{equation}
In other words $\mu_\lambda$ is the first node in the path from $\lambda$ to the root that is a right node or $\mu_\lambda = \mathcal{R}$ in case $\lambda = (j,0)$. Then we consider the sets 
\[\Upsilon(\lambda) = \{ \mu \in \Lambda_{I,0} \: : \: \lambda \succeq \mu \succeq \mu_\lambda\}. \]
Now let $\mathcal{T} \subset \Lambda_{I,0}$ be a tree of wavelet indices. One easily checks that for two siblings $\eta_1,\eta_2 \in \mathcal{T}$ it holds  $\Upsilon(\eta_1) \cap  \Upsilon(\eta_2) = \emptyset$ and $\bigcup_{\lambda \in \mathcal{V}(\mathcal{T})} \Upsilon(\lambda) = \mathcal{T}$. Consequently, all the conditions concerning $\Upsilon$ are fulfilled, see \eqref{eq:Upsilon}. In Figure \ref{Fig:hp-grid} a visualization for the sets $\Upsilon(\lambda)$ is provided.
\begin{figure}[!t]
	\centering
	\includegraphics[width=0.6\textwidth]{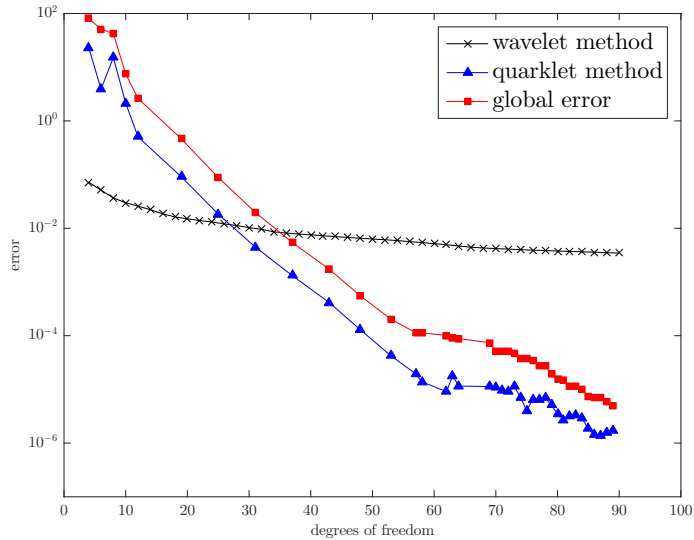}
	\caption{Error asymptotics for the test problem $f$ in semi-logarithmic scale. The black and blue lines depict the $L_2$-approximation error of the wavelet and quarklet method, respectively. The red line shows the behavior of the estimator $\mathcal{E}(T_N)^{1/2}$ in the quarklet case.}
	\label{Fig:conv_f}
\end{figure}
\begin{figure}[!t]
	\centering
	\begin{minipage}[b]{0.49\textwidth}
		\includegraphics[width=\textwidth]{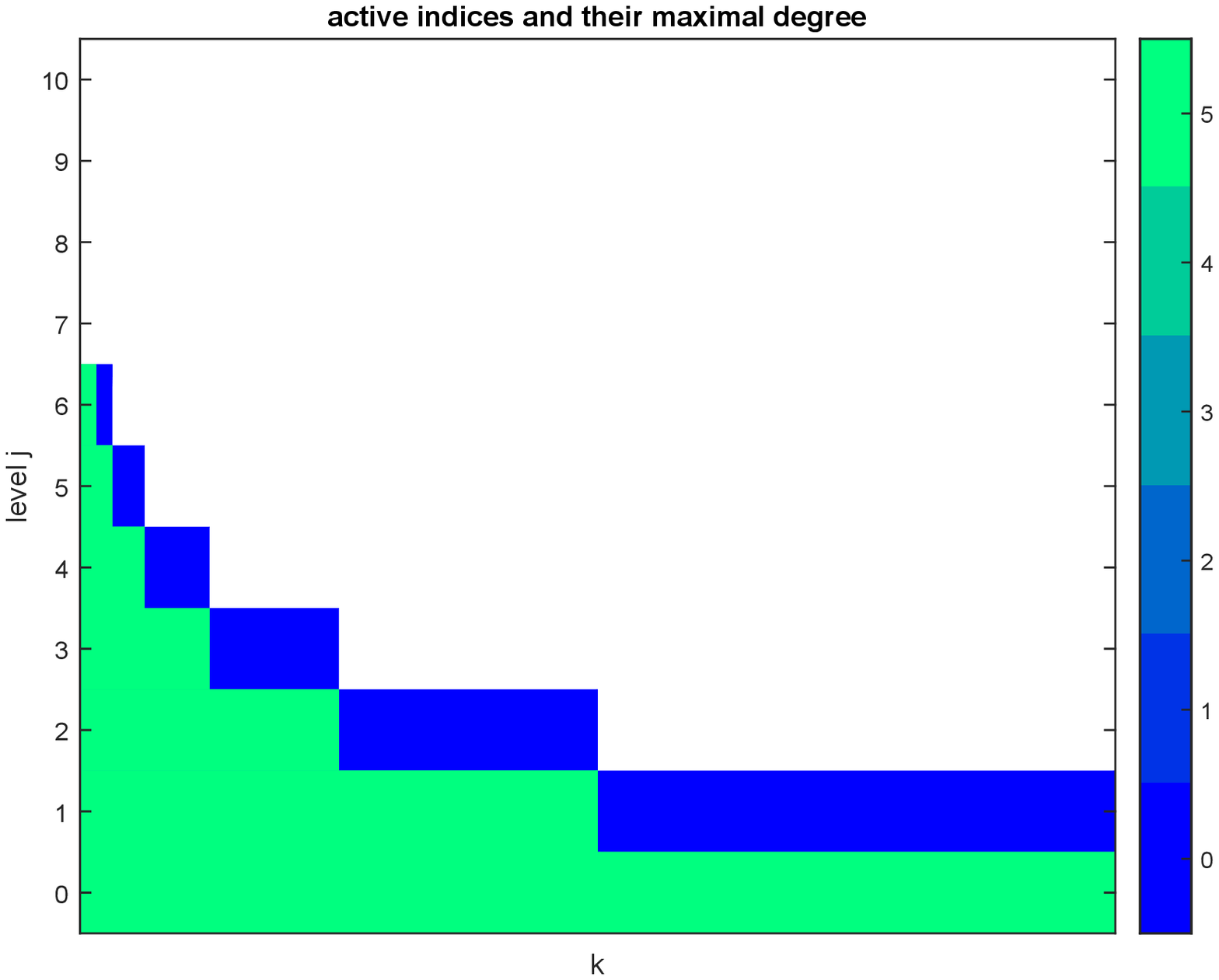}
	\end{minipage}
	\begin{minipage}[b]{0.49\textwidth}
		\includegraphics[width=\textwidth]{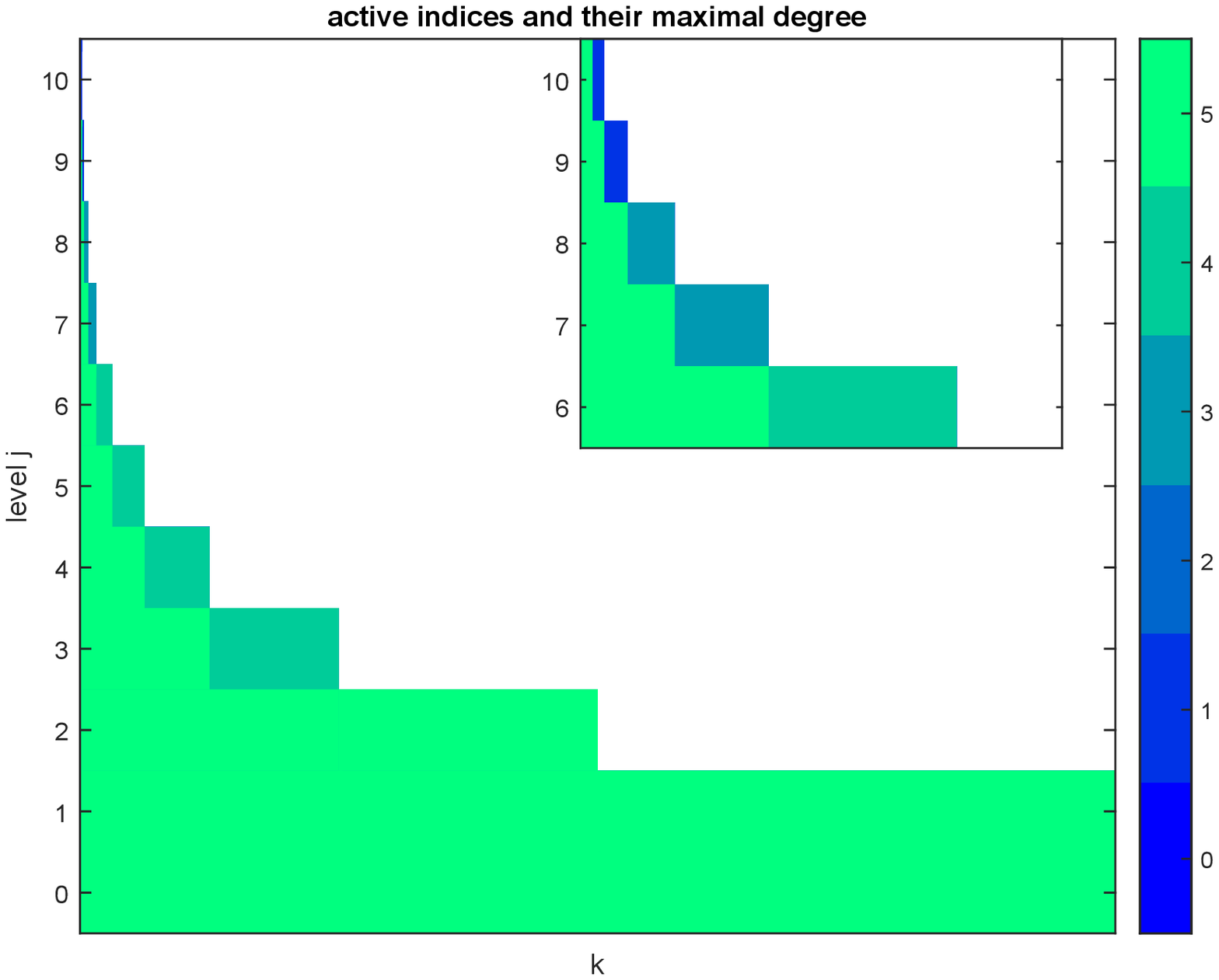}
	\end{minipage}
	\caption{Distribution of the active quarklet coefficients and their maximal polynomial degrees in the tree $T_N$ for the example $f$ after $N=12$ steps with $48$ degrees of freedom (left) and $N=50$ steps with $100$ degrees of freedom (right).}
	\label{Fig:coef_f}
\end{figure}
\begin{figure}[!t]
	\centering
	\includegraphics[width=0.6\textwidth]{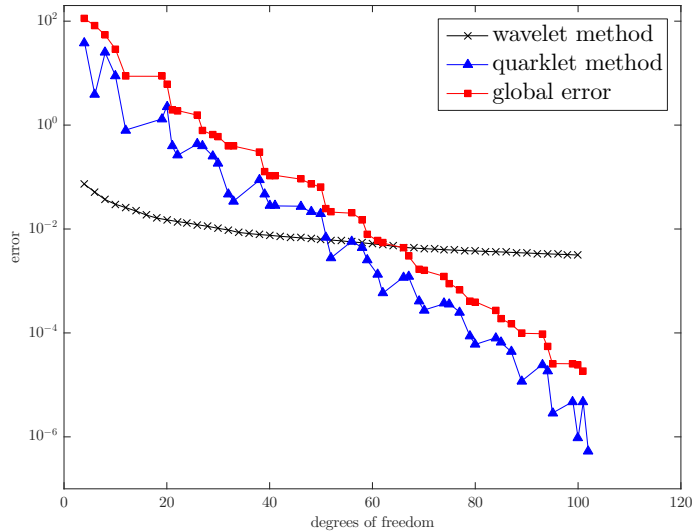}
	\caption{Error asymptotics for the test problem $g$ in semi-logarithmic scale. The black and blue lines depict the $L_2$-approximation error of the wavelet and quarklet method, respectively. The red line shows the behavior of the estimator $\mathcal{E}(T_N)^{1/2}$ in the quarklet case.}
	\label{Fig:conv_g}
\end{figure}
\begin{figure}[!t]
	\centering
	\includegraphics[width=0.49\textwidth]{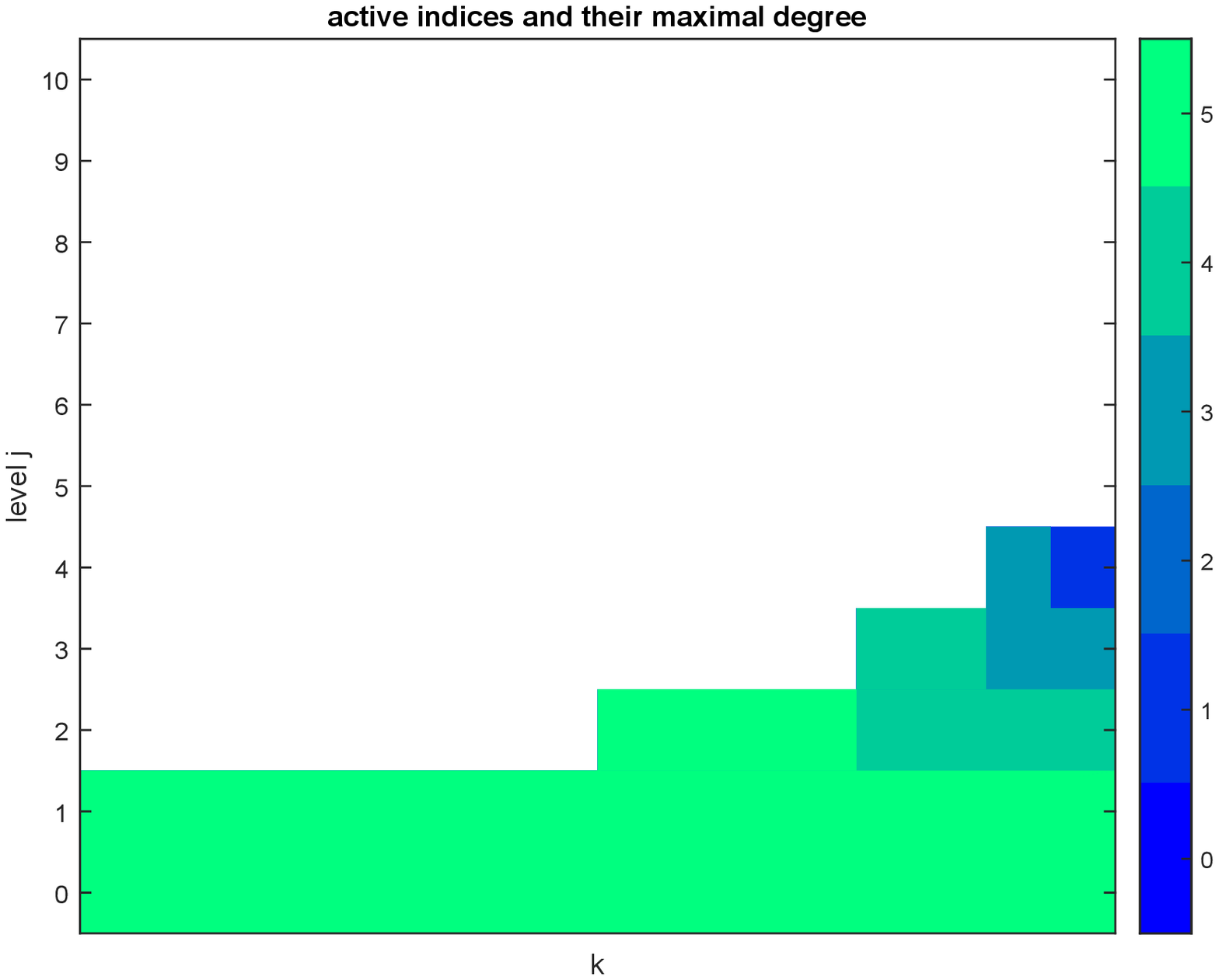}
	\includegraphics[width=0.49\textwidth]{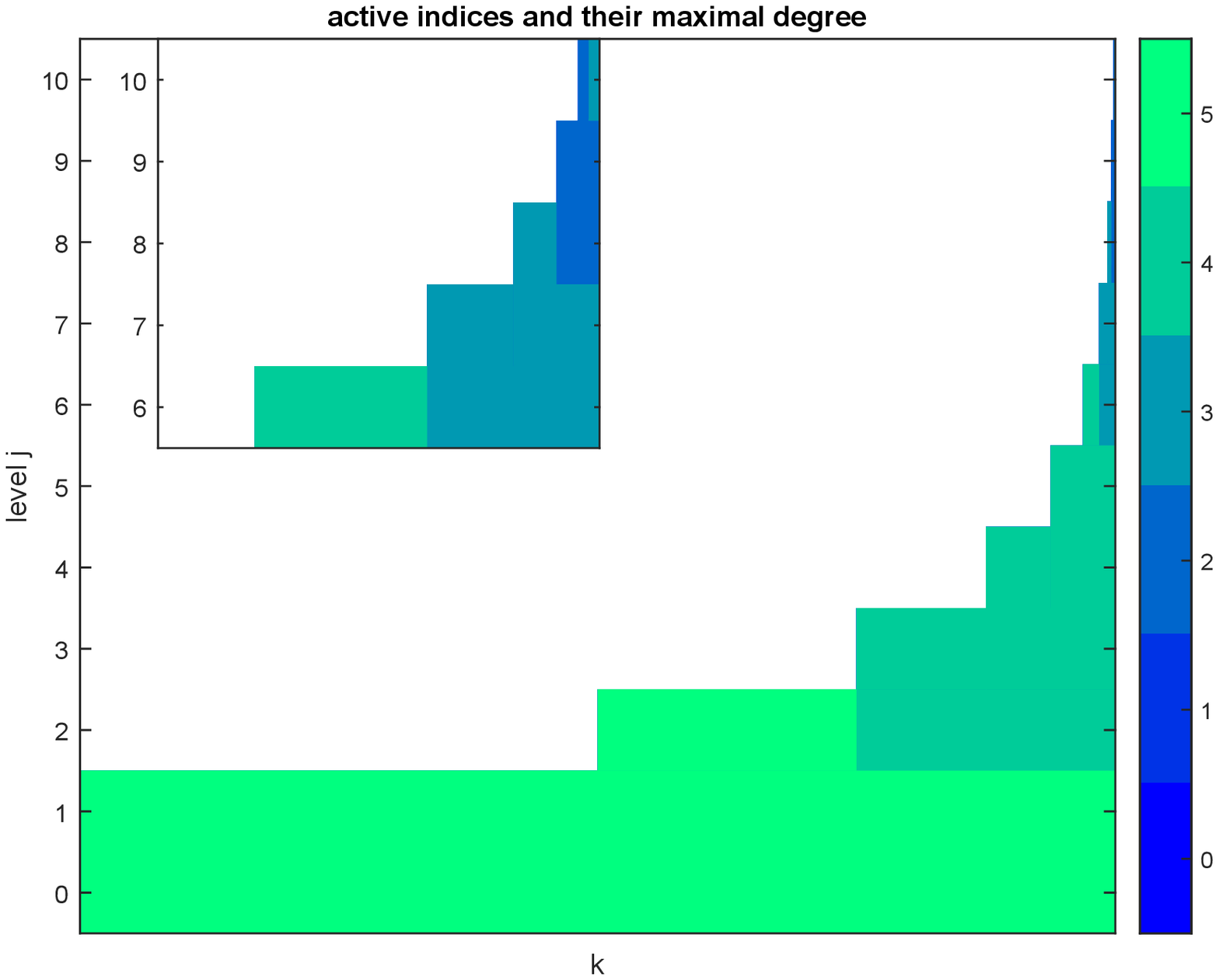}
	\caption{Distribution of the active quarklet coefficients and their maximal polynomial degrees in the tree $T_N$ for the example $g$ after $N=23$ steps with $50$ degrees of freedom (left) and $N=50$ steps with $102$ degrees of freedom (right).}
	\label{Fig:coef_g}
\end{figure}
\begin{figure}[!t]
	\centering
	\includegraphics[width=0.6\textwidth]{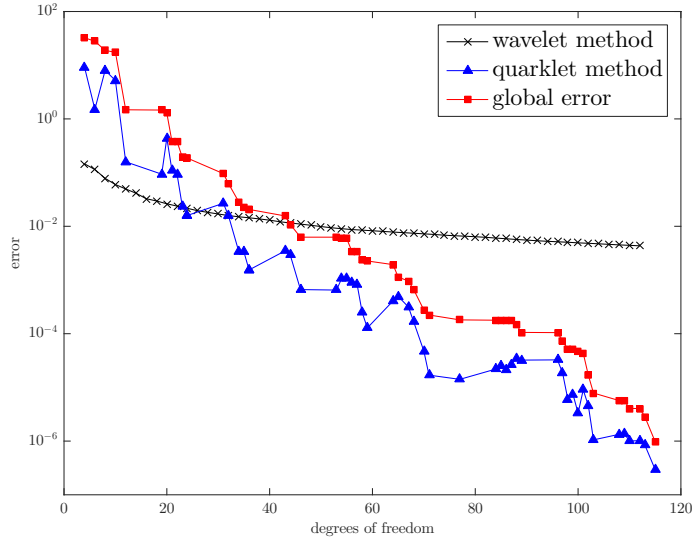}
	\caption{Error asymptotics for the test problem $u$ in semi-logarithmic scale. The black and blue lines depict the $L_2$-approximation error of the wavelet and quarklet method, respectively. The red line shows the behavior of the estimator $\mathcal{E}(T_N)^{1/2}$ in the quarklet case.}
	\label{Fig:conv_u}
\end{figure}
\begin{figure}[!t]
	\centering
	\includegraphics[width=0.49\textwidth]{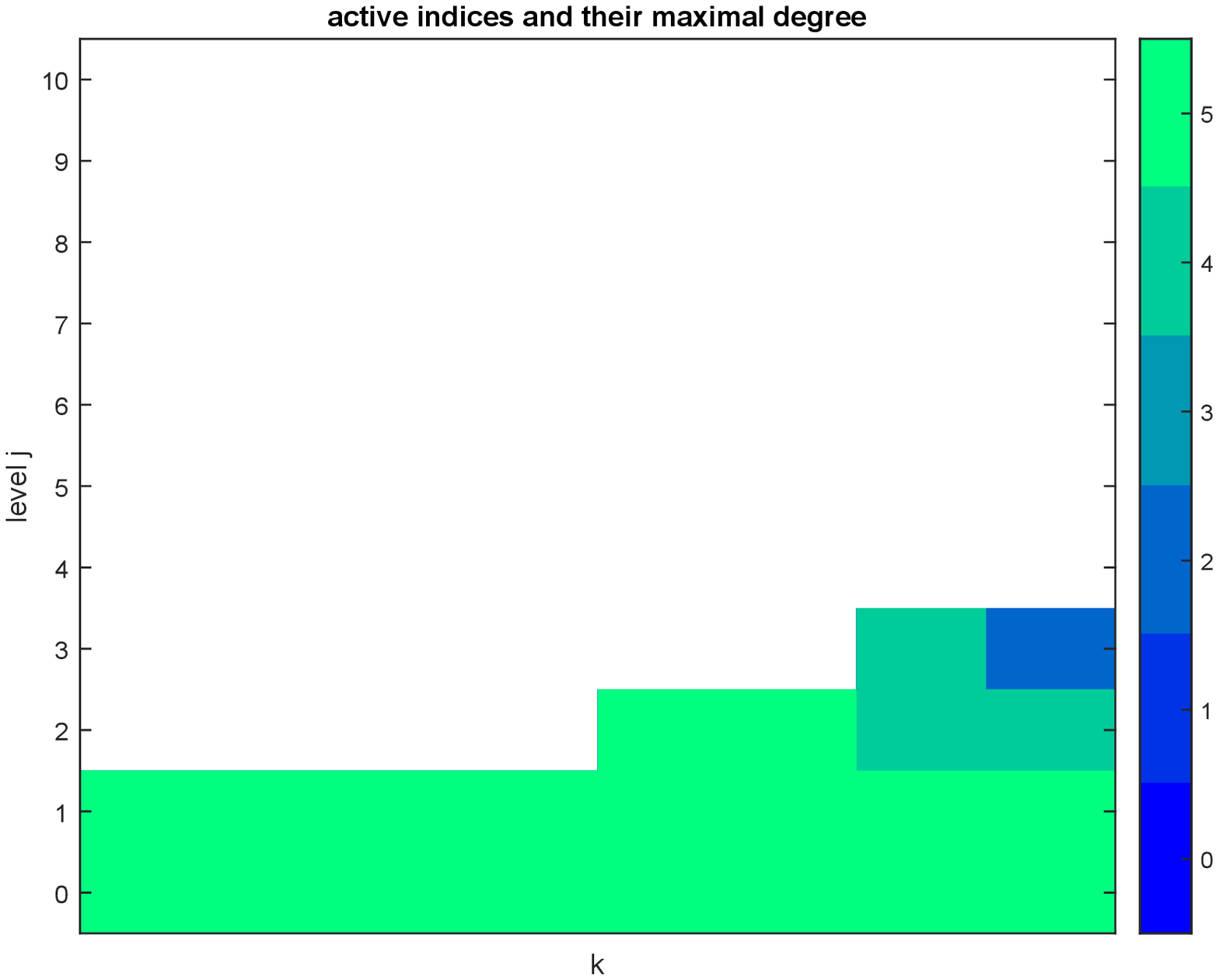}
	\includegraphics[width=0.49\textwidth]{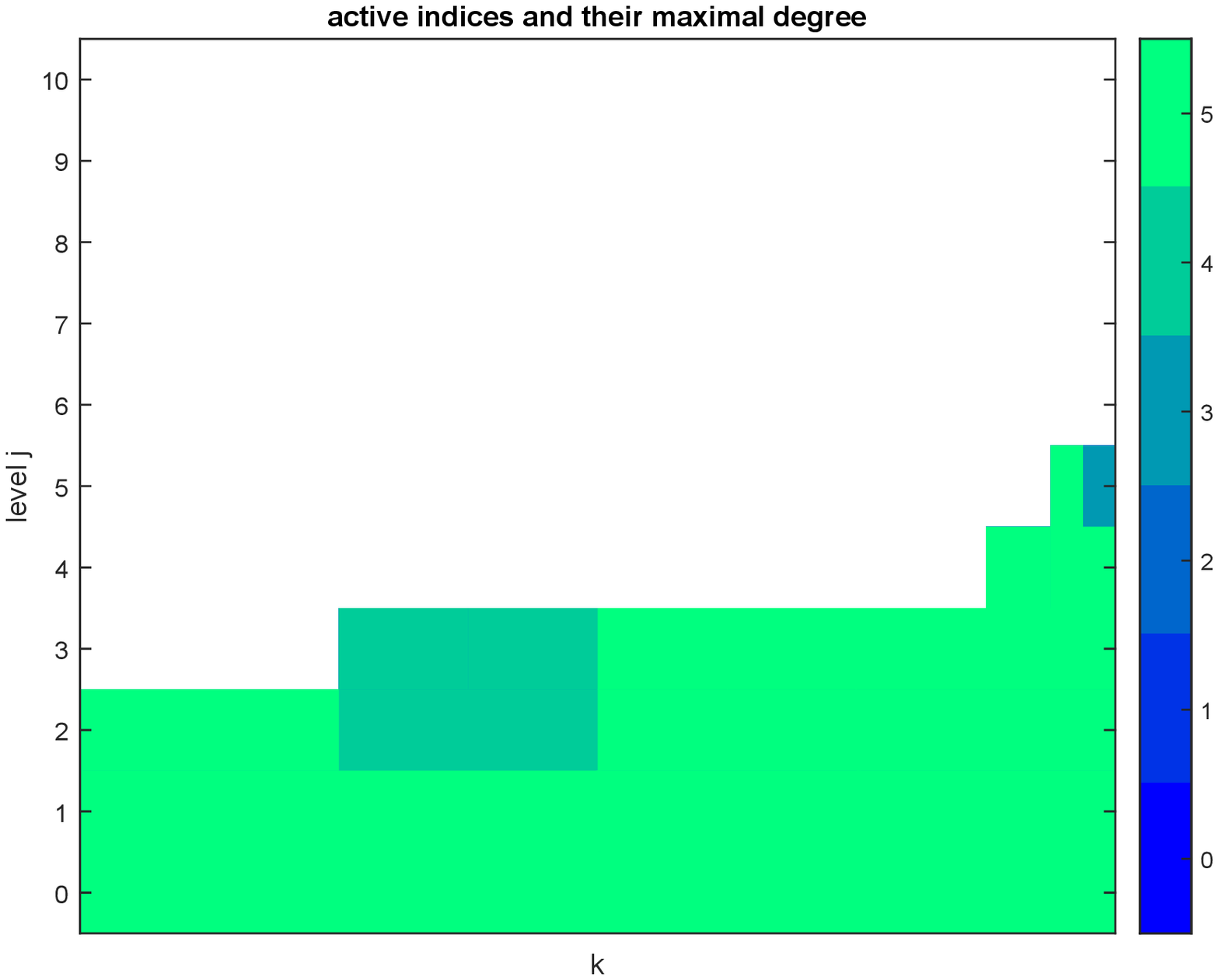}
	\caption{Distribution of the active quarklet coefficients and their maximal polynomial degrees in the tree $T_N$ for the example $u$ after $N=20$ steps with $43$ degrees of freedom (left) and $N=50$ steps with $103$ degrees of freedom (right).}
	\label{Fig:coef_u}
\end{figure} 
\begin{figure}[!t]
	\centering
	\includegraphics[width=0.6\textwidth]{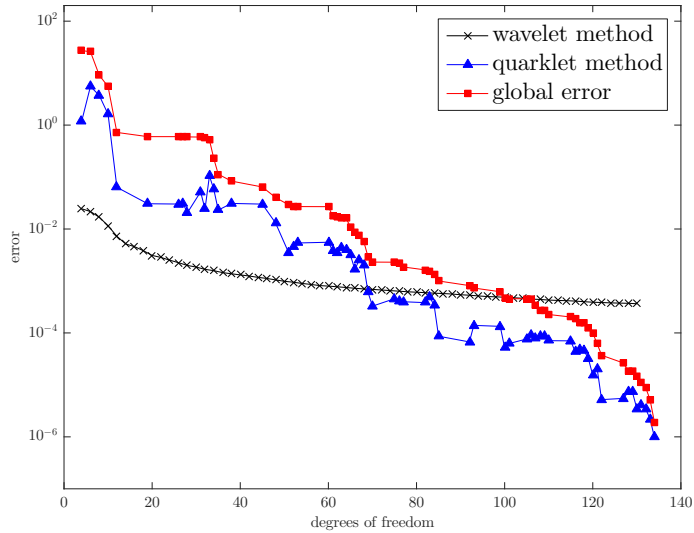}
	\caption{Error asymptotics for the test problem $v$ in semi-logarithmic scale. The black and blue lines depict the $L_2$-approximation error of the wavelet and quarklet method, respectively. The red line shows the behavior of the estimator $\mathcal{E}(T_N)^{1/2}$ in the quarklet case.}
	\label{Fig:conv_v}
\end{figure}
\begin{figure}[!ht]
	\centering
	\includegraphics[width=0.49\textwidth]{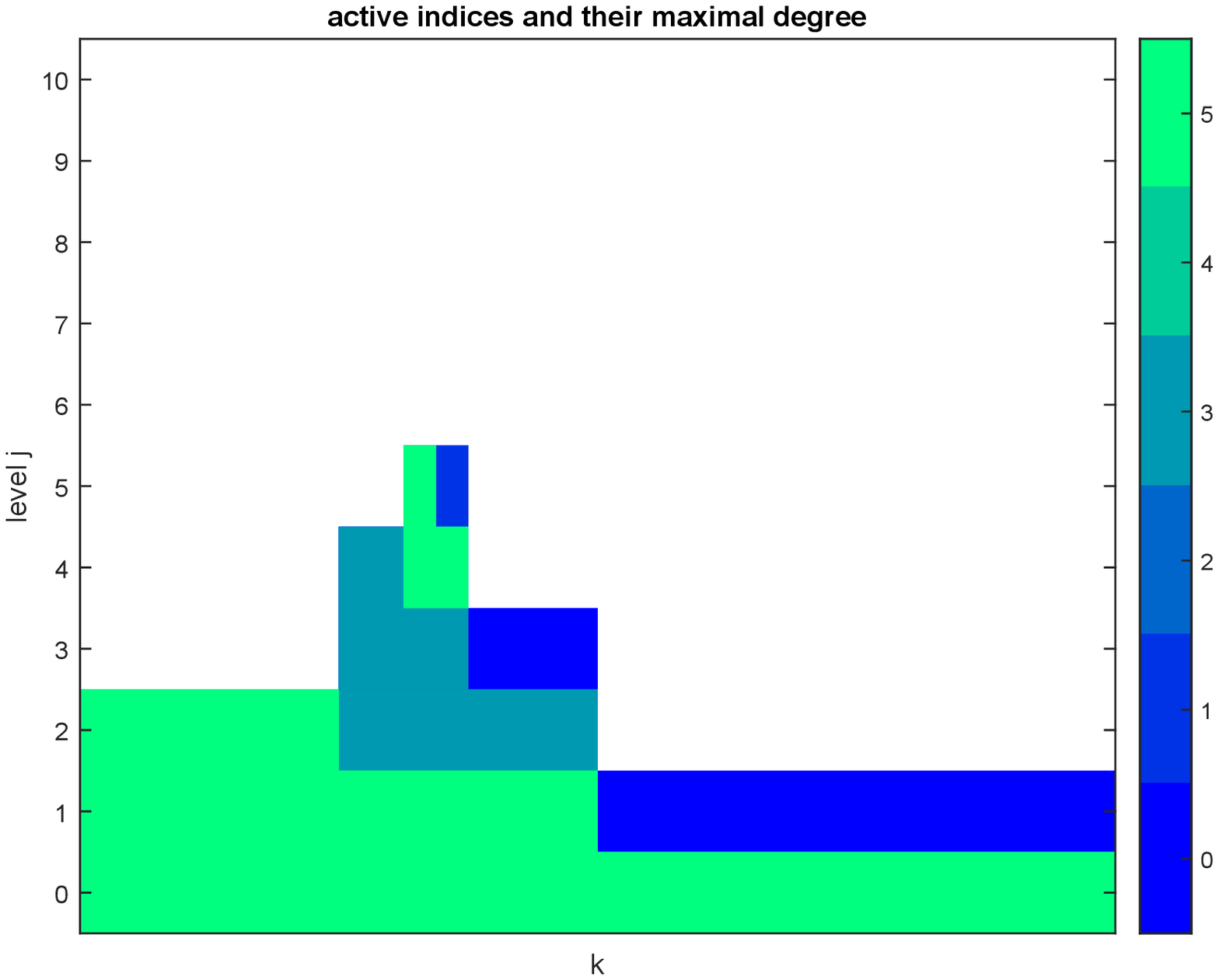}
	\includegraphics[width=0.49\textwidth]{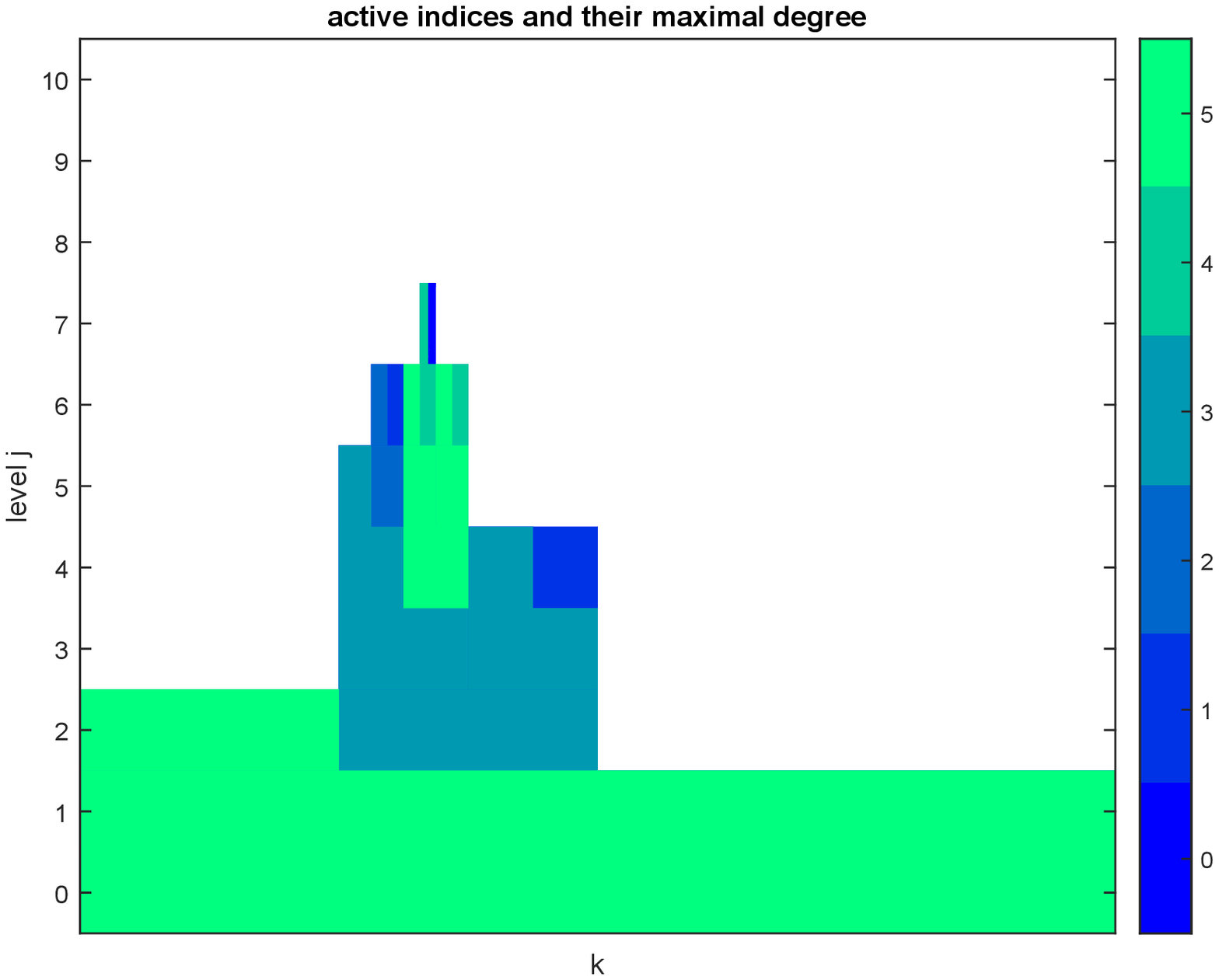}
	\caption{Distribution of the active quarklet coefficients and their maximal polynomial degrees in the tree $T_N$ for the example $v$ after $N=20$ steps with $52$ degrees of freedom (left) and $N=50$ steps with $110$ degrees of freedom (right).}
	\label{Fig:coef_v}
\end{figure}

We are now ready to apply our algorithm. The first example is the aforementioned function 
\[f(x) \coloneqq x^\alpha, \qquad \alpha > \frac12.\]
As usual we choose $\alpha = \frac{3}{4}$ for our investigations. Notice that this function is smooth except for $x = 0$. Therefore we expect many refinements in space towards the left boundary of the interval while the smooth part of $f$ should be well resolved by quarklets with a high polynomial degree. In Figure \ref{Fig:conv_f} one can observe the decay of the approximation error $\Vert f - f_{T_N} \Vert_{ L_2(I) }$ and the global error $\mathcal{E}(T_N)^{1/2}$ with increasing $N$ in semi-logarithmic scale. Recall that the degrees of freedom depend on $N$, see Lemma \ref{lemma:T_N_cardinality}.  We observe that the global error provides a very precise estimate of the $L_2$-approximation error, cf. Lemma \ref{lem_L2_globerr}. Moreover, in the semilogarithmic scale, the error decays linearly which indicates exponential convergence of type $e^{-\beta n^{\gamma}}$ with $\gamma$ close to $1$. As a comparison we also included the error of an adaptive wavelet method. To this end we  applied {\bf NEARBEST\textunderscore TREE} with the index set $\Lambda_I$ truncated at $p_{\textrm{MAX}}^\prime = 0$ and $j_{\textrm{MAX}} = 10$. This space adaptive scheme realizes the linear convergence rate $n^{-1}$. During the first steps of the algorithm the wavelet method achieves a higher accuracy. This implies that the choice of coefficients we used in the representation \eqref{rep_L2_gen} for the quarklet scheme can be improved. Nevertheless we observe the significantly higher asymptotic convergence rate of the adaptive quarklet method. In Figure \ref{Fig:coef_f} we present the active quarklet coefficients in $T_N$ at different stages of the algorithm {\bf NEARBEST\textunderscore TREE}. We notice a strong refinement in scale toward the singularity while also high polynomial degrees are used. In our second experiment we will investigate the role of the specific choice of the sets $\Upsilon(\lambda)$. To this end we will consider the reflected and shifted version of $f$ given by
\[g(x) \coloneqq (-x+1)^\alpha, \qquad \alpha > \frac12.\]
Again we set $\alpha = \frac34$. The results are depicted in Figures \ref{Fig:conv_g} and \ref{Fig:coef_g}. In comparison to the previous example we notice some slight improvements in the resulting convergence rate and the distribution of active coefficients. Especially we see that the higher polynomial degrees are more concentrated on the lower levels. In our previous example this is not the case. This is a consequence of the definition of the sets $\Upsilon(\lambda)$. A high polynomial degree on the root $\mathcal{R} = (0,0)$ automatically implies a high polynomial degree on the leftmost leaf and the nodes along this path. This is exactly the situation we encountered in Figure \ref{Fig:coef_f}. Another effect is the much faster allocation of degrees of freedom with increasing $N$ in the first steps of our first example. In contrast this occurs in a more uniform fashion for $g$. Indeed, it seems to be advantageous to have `smaller' sets $\Upsilon(\lambda)$ in areas where the function can be classified as `rough'. In this sense our choice of the sets $\Upsilon(\lambda)$ seems to be well suited for $g$. Nevertheless, we see that the exponential convergence shows up in both cases. In our final experiments we will inspect two more functions that are popular test cases as solutions to elliptic problems. The first is given by
\[u(x) \coloneqq 4 \frac{e^{ax}-1}{e^a-1}\left( 1 - \frac{e^{ax}-1}{e^a -1} \right). \]
As in \cite{bib:BBCCDDU} we choose $a=5$, although other values are possible. This function is smooth but has a large gradient at $x=1$. The second function is taken from \cite{Wih11} and exhibits a spike at $x= \frac13$. It is given by
\[v(x) \coloneqq \frac{x(1-x)}{1+10^4\left( x-\frac13\right) ^2}.\]
The resulting convergence rates and the distribution of active indices are depicted in Figures \ref{Fig:conv_u}-\ref{Fig:coef_v}. In particular, we notice that exponential decay of the errors is also achieved in these examples.

\vspace{0,3 cm}

\textbf{Funding.} This paper is a result of the DFG project `Adaptive high-order quarklet frame methods for elliptic operator equations' with grant numbers DA360/24$-$1 (Stephan Dahlke, Marc Hovemann) and RA2090/3$-$1 (Thorsten Raasch).

\end{document}